\documentclass[12pt]{amsart}

\usepackage{amsmath,amsfonts,amssymb,amscd}
\begin{document}
\renewcommand{\thefootnote}{\fnsymbol{footnote}}
\pagestyle{plain}

\title{A stronger concept of K-stability}
\author{Toshiki Mabuchi${}^*$}
\maketitle
\abstract
In this paper, by introducing a wider class of one-parameter group actions
for test configurations, we have a stronger version of the definition of
K-stability.
This allows us to obtain some key step of \cite{M4} in proving
that constant scalar curvature polarization implies K-stability 
for polarized algebraic manifolds.
\endabstract

\section{Introduction}

In this paper, we fix once for all a polarized algebraic manifold $(M,L)$
consisting of a connected projective algebraic manifold $M$ of dimension $n$, defined over $\Bbb C$, 
and a very 
ample line bundle $L$ over $M$. Take a 1-dimensional algebraic torus 
$T_0 :=  \Bbb C^*$ which acts on the affine space
$$
\Bbb A^1\; = \,\{\,z \in \Bbb C\,\}
$$
by multiplication of complex numbers with some positive weight $\alpha$, so that the action 
is given by${}^{\dagger}$
\footnotetext{ ${}^{\dagger}$ 
In Donaldson's definition of a test configuration, 
$(\mathcal{M}_z, \mathcal{L}_z)  \cong \; (M,L^{\alpha})$, $z \neq 0$. 
for some exponent $\alpha$, where the weight is chosen to be 1. 
However for our definition, the weight is $\alpha$ and the exponent is 1. 
When $T_0$ is replaced by its unramified cover of degree $\alpha$ in Donaldson's definition, 
the pair of the weight and the exponent changes from $(1, \alpha )$ to $(\alpha , 1)$, yielding our definition equivalent to Donaldson's.}
$$
\Bbb C^* \,(= T_0)\, \times\,  \Bbb A^1\, \to\, \Bbb A^1, \;\quad (t, z) \to t^{\alpha} z.
$$
Let $(\mathcal{M}, \mathcal{L})$ be a test configuration
(cf. Donaldson \cite{D1}) 
for $(M,L)$, so that we have a $T_0$-equivariant projective morphism 
$$
\pi : \mathcal{M} \to \Bbb A^1\; (= \,\{\,z \in \Bbb C\,\}\,)
$$ 
of an irreducible reduced algebraic variety $\mathcal{M}$, defined over $\Bbb C$,
onto the affine space $\Bbb A^1$,
where $\mathcal{L}$ is a relatively very ample
 invertible sheaf on the fiber space $\mathcal{M}$ over $\Bbb A^1$. 
Then the restriction  of the pair
$(\mathcal{M}, \mathcal{L})$ to each fiber $\mathcal{M}_z := \pi^{-1}(z)$ admits
a holomorphic isomorphism
$$
(\mathcal{M}_z, \mathcal{L}_z) \; \cong \; (M,L),
\qquad z\neq 0,
$$
where the $T_0$-action on $\mathcal{M}$
lifts to a $T_0$-linearization of 
$\mathcal{L}$.
We now consider the identity component
$$
\mathcal{P}:= \operatorname{Aut}^0(\mathcal{M})^{T_0}
$$
of the group of all holomorphic automorphisms of $\mathcal{M}$ commuting with the $T_0$-action on $\mathcal{M}$. Note that every element of $\mathcal{P}$ maps fibers of 
$\pi$ to fibers of $\pi$, inducing  multiplication by a complex number
on the base space $\Bbb A^1$. 
For the Lie algebra $\frak p$ of $\mathcal{P}$,
in view of the $T$-equivariance of $\pi$,
we can choose an element $X_0\neq 0$ in $\frak p$ 
generating holomorphically the $T_0$-action on $\mathcal{M}$ such that
$\exp \, (2\pi \sqrt{-1}  \, X_0 /\alpha ) \,=\, \operatorname{id}_{\mathcal{M}}$ and that
$$
\pi_* X_0 \; =\; \alpha\,z\,\frac{\partial}{\partial z},
$$
where $T_0$ is often identified with 
$\{\,\exp (s X_0)\,;\, s\in \Bbb C\,\}$.
For the fiber space $\mathcal{M}$ over $\Bbb A^1$, 
consider the subgroup $\operatorname{Aut}(\mathcal{M}, \Bbb A^1)^{T_0}$ in $\mathcal{P}$ of
all fiber-preserving elements in  $\mathcal{P}$.
Let $\mathcal{Q} = \mathcal{Q}(\mathcal{M})$ 
be its identity component
$$
\mathcal{Q}:= \operatorname{Aut}^0(\mathcal{M}, \Bbb A^1)^{T_0}.
$$
Then by restricting each automorphism in  $\mathcal{Q}$ to the fiber 
$\mathcal{M}_1  \, (\, =\,\pi^{-1}(1)$ $=\, M)$, 
$\mathcal{Q}$ is viewed as a closed algebraic subgroup 
$$
\mathcal{Q} \; \hookrightarrow \; \operatorname{Aut}^0(M) \, (=\, \operatorname{Aut}^0(\mathcal{M}_1)),
$$
where  injectivity follows from the commutativity of the elements in $\mathcal{Q}$ 
with the $T_0$-action on $\mathcal{M}$. 
Define $H = H(\mathcal{M})$ as the intersection
$$
H\; :=\; \mathcal{Q}\,\cap \, G, 
$$
where $G$ is the maximal connected linear algebraic subgroup of $\operatorname{Aut}^0(M)$.
Let  $\frak g$ and
$\frak h \,=\, \frak h (\mathcal{M})$ 
be the Lie algebras of $G$ and $H = H(\mathcal{M})$, respectively.
By fixing a K\"ahler form $\omega$ in the class $c_1(L)_{\Bbb R}$,
we say that an element $Y$ of $\frak g$ is {\it Hamiltonian} if
$$
i (Y)\, \omega \; =\; \bar{\partial} f_{\omega, Y}
\leqno{(1.1)}
$$
for some real-valued smooth function $f_{\omega,Y} \in C^{\infty}(M)_{\Bbb R}$
with normalization condition $\int_M f_{\omega,Y}\, \omega^n = 0$.
Then for every Hamiltonian element $Y$ of $\frak g$, the closure 
in $G$ of the holomorphic one-parameter 
group generated by $Y$
is an algebraic torus. 
In view of the inclusion $\frak h \subset \frak g$, 
we define a subset $\,\frak S\, =\, \frak S (\mathcal{M})\,$ of $\frak h$ as
$$
\frak S  \,:\; \text{the set of all Hamiltonian elements in $\frak h$}. 
$$
Since the $H$-action on 
$\mathcal{M}$ lifts to an $H$-linearization of 
some positive integral multiple of $\mathcal{L}$, there exists a finite  
unramified cover $\tilde{H}$
of $H$
such that the $\tilde{H}$-action on $\mathcal{M}$ 
induced by 
$$
H\; \subset \; \mathcal{Q} \; =\;
 \operatorname{Aut}^0(\mathcal{M}, \Bbb A^1)^{T_0}
$$
lifts to an $\tilde{H}$-linearization of $\mathcal{L}$.
Given a element $Y$ of $\frak S$, the closure 
$\mathcal{S}:=\bar{\tau}_Y$ in $\tilde{H}$ of the holomorphic one-parameter 
group 
$$
\tau_Y\, :=\,\{\, \exp (t Y)\in \tilde{H}\,;\, t\in \Bbb C\,\}
$$
is an algebraic torus. Put
 $\delta_Y := \dim_{\Bbb C}\,\bar{\tau}_Y$.
If $\delta_Y \leq 1$, 
then $\bar{\tau}_Y$ and ${\tau}_Y$ coincide, and  
in this case, we say that $Y$ is {\it quasi-regular}.
On the other hand, if $\delta_Y > 1$, then $Y$ is  
said to be {\it irregular}.
These concepts are analogous to the quasi-regularity or irregularity 
in Sasakian geometry.

\medskip
Fix an element $Y$ of $\frak S$. Then $X:= X_0 + Y \in \frak p$ 
generates a holomorphic one-parameter group 
$$
T\; := \; \{\,\exp (tX)\,;\, t\in \Bbb C\,\}
$$
in ${\Sigma}:= {\mathcal{S}}\times T_0$ by 
$\exp (tX) = \exp (tY)\cdot\exp (tX_0)$.
Here $\exp (tY)$ is regarded as an element of  ${\mathcal{S}}$,
while $\exp (tX_0)$ sits in $T_0$.
Since the $T_0$-action on $\mathcal{M}$ lifts to a $T_0$-linearization of 
$\mathcal{L}$, 
a similar lifting is true also for the $T$-action.
Hence, whether $Y$ is quasi-regular or irregular,
${\Sigma}$  and its subgroup $T$
act on the vector bundles${}^{\dagger}$ 
\footnotetext{ ${}^{\dagger}$  
Let $h_{m;1}$ be a Hermitian metric on the 
vector space $(E_m)_z$ at $z =1$.
By \cite{D3} and \cite{Mas}, 
we have a ${\Sigma}$-equivariant (and hence $T$-equivariant) 
trivialization for $E_m$,
$$
E_m \; \cong \; \Bbb A^1 \times (E_m)_0,
\leqno{(1.2)}
$$
taking $h_{m;1}$ to a Hermitian metric $h_{m;0}$ on $(E_m)_0$
which is preserved by the action of the compact torus $(T_0)_c \cong S^1$ in $ T_0$.
For the maximal compact subgroup $\mathcal{S}_c$ 
of $\mathcal{S}$, by taking an average, we may also assume that 
both $h_{m;0}$ and $h_{m;1}$ are $\mathcal{S}_c$-invariant.} 
$$ 
E_m:= \pi_*\mathcal{L}^m, \qquad m =1,2,\dots,
$$
over $\mathcal{M}$. 
Here vector bundles and invertible sheaves 
are used interchangeably throughout this paper.
 Put $N_m := \dim H^0 (M, \mathcal{O}^{}_M(L^m))$. 
Replacing $\mathcal{L}$ by its suitable positive integral multiple 
if necessary,
we may assume that 
$\dim H^0(\mathcal{M}_0, \mathcal{L}_0^m)$ is $N_m$ for all $m$ and that
$$
\otimes^m H^0(\mathcal{M}_z, \mathcal{L}_z) 
\to H^0(\mathcal{M}_z, \mathcal{L}_z^m), 
\qquad m = 1,2,\dots, 
$$
are surjective for all $z\in \Bbb A^1$ (cf. \cite{M1.5}, Remark 4.6).
In view of
$$
\pi_* X \; =\;  \alpha\, z\,\frac{\partial}{\partial z},
$$
the projective morphism $\pi : \mathcal{M}\to \Bbb A^1$ is  
always $T$-equivariant, though $T$ is not isomorphic to $\Bbb C^*$ 
in irregular cases.
Then a triple $(\mathcal{M},\mathcal{L},X)$ is called a {\it 
generalized test configuration of $(M,L)$}, 
if $Y:= X - X_0$ belongs to $\frak S = \frak S (\mathcal{M})$.
Note that a test configuration (cf. \cite{D1}) of $(M, L)$ in an ordinary sense corresponds to 
the triple $(\mathcal{M},\mathcal{L},X_0)$, i.e., to the case 
$Y =0$. 
Now for the triple 
$(\mathcal{M},\mathcal{L},X)$, we consider the weight 
$$
w_m = w_m(\mathcal{M},\mathcal{L},X)
$$
of the $T$-action on $\det (E_m)_0$, i.e., the trace of the endomorphism on 
the vector space $(E_m)_0$ 
induced by $X$.
Now for $m \gg 1$, 
we have the expansion
$$
\frac{w_m}{m N_m} \; = \; F_0 + F_1 m^{-1} + F_2 m^{-2} + \dots ,
$$
with coefficients $F_i = F_i (\mathcal{M},\mathcal{L},X)\in \Bbb R$.
We now say that $(M,L)$ is {\it K-semistable in a strong sense}
if  
$$
F_1(\mathcal{M},\mathcal{L},X) \; \leq \; 0
$$
for all generalized test configurations $(\mathcal{M},\mathcal{L},X)$ 
of $(M,L)$. 
Moreover, a K-semistable $(M,L)$ is said to be
{\it K-stable in a strong sense} if, for every generalized test configurations 
$(\mathcal{M},\mathcal{L},X)$ of $(M,L)$,
$$
\text{$\mathcal{M} = \Bbb A^1\times M$ if and only if 
$F_1(\mathcal{M},\mathcal{L},X)$ vanishes,}
$$
where $T$ does not necessarily act on the second factor $M$ of $\Bbb A^1\times M$
trivially. Note that K-stability in a strong sense implies
K-stability in an ordinary sense, 
and vice versa (cf. \cite{MN}). 
 In this paper, by effectively using the concept of generalized test configurations, 
 we shall complete the proof of the following:

\medskip\noindent
{\bf Main Theorem.}
{\em A polarized algebraic manifold $(M,L)$ is K-stable if
the class $c_1(L)_{\Bbb R}$ admits a K\"ahler metric  
of constant scalar curvature.} 

\medskip
In Case 1 of Section 5  of \cite{M4}, some key step in the proof of Main Theorem 
above is not given. We here give the omitted proof by showing that, 
given a K\"ahler metric in $c_1(L)_{\Bbb R}$ of constant scalar curvature, 
the vanishing of $F_1(\mathcal{M}, \mathcal{L})$ yields
a suitable generalized test configuration $(\mathcal{M}', \mathcal{L}', X')$,
from which $\mathcal{M} = \mathcal{M}'\cong \Bbb A^1\times M$ easily follows.

\medskip
\section{A generalized test configuration associated to $W$}

In this section, we fix a test configuration 
$(\mathcal{M}, \mathcal{L})$ 
 of $(M, L)$. 
To each Hamiltonian element $W$ in $\frak g$,
 we associate a generalized test configuration $(\mathcal{M}', \mathcal{L}', X')$ 
of $(M, L)$ as follows:

\medskip\noindent
Since the $T_0$-action on $E_m$ preserves the fiber 
$(E_m)_0$ of $E_m$ over the origin, we write the associated 
representation on $(E_m)_0$ as
$$
\Psi^{}_{m, X_0} (=\Psi^{}_{m, X_0;\mathcal{M}, \mathcal{L}}) : 
T_0 \to \operatorname{GL}((E_m)_0).
$$
For each $\hat{t}\in T_0 $, we set $t:= \hat{t}^{N_m}\in T_0$,
and define an algebraic group homomorphism 
${\Psi}^{\operatorname{SL}}_{m,X_0}\,
(= {\Psi}^{\operatorname{SL}}_{m,X_0;\mathcal{M}, \mathcal{L}}) : 
T_0  \to \operatorname{SL}((E_m)_0)$ by 
$$
{\Psi}^{\operatorname{SL}}_{m, X_0}(\hat{t})\;  :=\; 
\frac{\Psi^{}_{m, X_0} (t)}
 {\det (\Psi^{}_{m, X_0} (\hat{t}))},  
\qquad \hat{t}\in T_0.
$$
Let $\Pi_m : \operatorname{GL}((E_m)_0) \to \operatorname{PGL}((E^*_m)_0)$
be the natural projection induced by the contragradient representation.
For each $s \in \Bbb C$, put $\hat{s}:= s/N_m$ and $z:= e^s$.  
In view of the identification of $T_0$ with $\Bbb C^*$, by setting
$$
\mu^{}_{m,s} \;:= \; (\Pi^{}_m \circ 
{\Psi}^{\operatorname{SL}}_{m, X_0}) (z),\,
$$
we define 
$M_{s} :=  \mu^{}_{m,s} (M)$, 
$s \in \Bbb C$, where $M$
is regarded as a submanifold of $\Bbb P^*((E_m)_0)$ 
by the inclusion
$$
M\,=\, \mathcal{M}_1 \subset 
\Bbb P^*((E_m)_1) \cong \Bbb P^*((E_m)_0).
\leqno{(2.1)}
$$ 
Here the last isomorphism is
induced by the 
$T_0$-equivariant trivialization 
$E_m \cong \Bbb A^1\times (E_m)_0$ in (1.2).
Hence identifying the projective bundle $\Bbb P^*(E_m)$ with the product 
$\Bbb A^1 \times \Bbb P^*((E_m)_0)$,
we obtain
$$
\mathcal{M}_z \; =\;  \{z\}\times M_s, \qquad s \in \Bbb C,
\leqno{(2.2)}
$$
where $\mathcal{M}$ is regarded as a subvariety of $\Bbb P^*(E_m)$.
Since $W$ is a Hamiltonian element of $\frak g$, we see that 
$$
Y':= - W
$$ 
is also Hamiltonian.
For the holomorphic one-parameter group $\tau^{}_{Y'}$ generated by $Y'$ in $G$, 
the closure $\mathcal{S}' := \bar{\tau}^{}_{Y'}$
in $G$ (later on, $\mathcal{S}'$ will be replaced by its finite unramified cover) 
is an algebraic torus. Put 
$$
\begin{cases}
&\hat{\kappa}_z \, := \; \Psi_{1,X_0} (z),\\
&\kappa_z \, := \; \mu_{1,s}\; =\;\Pi_1\circ\Psi_{1,X_0} (z),
\end{cases}
$$
for all $z \; (= e^s )\in \Bbb C^*$, where $m=1$.
Now for $m =1$ and $z = e^s$,
we define a subset $U$ of $(\Bbb A^1\setminus \{0\})\times \Bbb P^*((E_m)_0)$ 
by
$$
U:= \bigcup_{0\neq z\in \Bbb A^1} \; 
 \{z\} \times M_s\; \,
 \biggl ( =\; \bigcup_{s\in \Bbb C} \; 
 \{z\} \times M_s  \biggl )\; 
= \; \mathcal{M}\setminus \mathcal{M}_0.
$$
The algebraic group $\Sigma' :=\mathcal{S}'\times \Bbb C^* \, 
(= \mathcal{S}'\times T_0)$ 
acts biregularly on $U$ by
$$
\Sigma' \times U \to U,
\;\quad ((\theta , t), (z,p) ) \mapsto (tz, {\kappa}^{}_{tz}\theta \,{\kappa}_{z}^{-1}p),
\leqno{(2.3)}
$$
where $(\theta ,t)\in \mathcal{S}' \times \Bbb C^*\, (=\Sigma' )$ and $(z,p)\in U \subset (\Bbb A^1\setminus \{0\})\times \Bbb P^*((E_1)_0)$.
For each $\theta\in\mathcal{S}'$,
let $\hat{\theta}$ denote the element of $\operatorname{GL}((E_1^*)^{}_0) $ induced by 
$\theta$ via the identification  
$$
(E_1^*)^{}_0 \; \cong \; (E_1^*)^{}_1 \; \cong  \; H^0(M, \mathcal{O}^{}_M(L))^*
$$
obtained from (1.2) applied to $m=1$.
Then (2.3)  is induced by the $\Sigma'$-action on $(\Bbb A^1 \setminus \{0\})\times (E^*_1)_0 $ defined by
$$
((\theta,t),(z,q))\, \mapsto\, (tz, \hat{\kappa}^{}_{tz}\,\hat{\theta}\,\hat{\kappa}_{z}^{-1}q),
\leqno{(2.4)}
$$
where $(\theta,t)\in \mathcal{S}' \times \Bbb C^*\, (=\Sigma' )$ and $(z,q)\in (\Bbb A^1 \setminus \{0\})\times (E^*_1)^{}_0 $.
Since the line bundle $\mathcal{O}_{\Bbb P^* (E_1)}(-1)$ 
is viewed as the blowing-up of $\Bbb A^1 \times (E^*_1)_0 $ along 
$\Bbb A^1 \times \{0\}$,
we see that
 the $\Sigma'$-action in (2.3)  naturally lifts to 
a $\Sigma'$-linearization of 
$\mathcal{L}\; =\; \mathcal{O}_{\Bbb P^* (E_1)}(1)_{|U}$
 over $U=\mathcal{M}\setminus \mathcal{M}_0$ induced by (2.4). 
 Let $\operatorname{pr}_1: (\Bbb A^1\setminus \{0\})\times \Bbb P^*((E_{1})_0)
\to \Bbb A^1\setminus \{0\}$ and 
$\operatorname{pr}_2: (\Bbb A^1\setminus \{0\})\times \Bbb P^*((E_{1})_0)
\to \Bbb P^*((E_{1})_0)$ be the projections to the first factor and to the second factor, 
respectively.
Take a $\Sigma'$-equivariant compactification $\bar{U}$ (cf.~\cite{S2}) of the 
algebraic variety $U$
such that the restriction of $\operatorname{pr}_1$ to $U$ extends to 
a $\Sigma'$-equivariant morphism
$$
\bar{\operatorname{pr}}: \bar{U} \to \Bbb P^1 (\Bbb C )
$$ 
onto $\Bbb P^1 (\Bbb C ):=\Bbb A^1 \cup \{\infty \}$. 
Here the first factor $\mathcal{S}'$ of $\Sigma'\,( = \mathcal{S}'\times T_0)$ 
acts on $\Bbb P^1 (\Bbb C )$ trivially.
By restricting $\bar{\operatorname{pr}}$ to 
$\hat{U}:= \bar{U} \setminus \bar{\operatorname{pr}}^{-1}(\{\infty \})$,
we obtain a $\Sigma'$-equivariant proper morphism
$$
\hat{\operatorname{pr}}: \hat{U} \to \Bbb A^1.
$$
Replacing $\hat{U}$ by its suitable $\Sigma'$-equivariant desingularization,
we may assume without loss of generality that $\hat{U}$ is smooth and sits over 
$\mathcal{M}$
by a $T_0$-equivariant modification
$$
\iota\, : \;\hat{U} \to \mathcal{M}.
$$
Fix a general hyperplane $\mathcal{H}$ on $\Bbb P^*((E_{1})_0)$. 
We also fix a divisor $\hat{\mathcal{H}}$ on $\hat{U}$ such that
$\hat{\mathcal{H}}
 - {\mathcal{H}} '$ is effective and that
$$
 \operatorname{Supp}( \hat{\mathcal{H}}
 - {\mathcal{H}} ')\; \subset \;\hat{\operatorname{pr}}^{-1}(0),
$$
where ${\mathcal{H}}'$
denotes the irreducible reduced divisor on $\hat{U}$ obtained as the
closure of  ${\operatorname{pr}}_2^{-1}(\mathcal{H})\cap U$ in $\hat{U}$.
Then
the divisor $\hat{\mathcal{H}}$ defines a line bundle 
$\hat{\mathcal{L}}:= \mathcal{O}_{\hat{U}}(\hat{\mathcal{H}})$ over $\hat{U}$ extending 
$\mathcal{L}_{|U}$. 
In view of \cite{Mu}, p.35,  
the $\Sigma'$-action on $\hat{U}$ lifts to
a $\Sigma'$-linearization of some positive integral multiple $\hat{\mathcal{L}}^{\beta}$ 
of $\hat{\mathcal{L}}$.
We then consider the unramified cover
$$
\Sigma' \to \Sigma', \qquad \sigma \mapsto \sigma^{\beta}.
$$
Given $\Sigma'$ on the right-hand side,
replacing this by $\Sigma'$ on the left-hand side (and hence 
the weight $\alpha$ in the introduction is replaced by $\alpha \beta$, 
so that the exponent can be chosen to be 1), 
we may assume from the beginning that the $\Sigma'$-action on $\hat{U}$ lifts to
a $\Sigma'$-linearization of $\hat{\mathcal{L}}$. 
For each $z \in \Bbb A^1\setminus \{0\}\subset \Bbb P^1(\Bbb C )$, 
let $\hat{\mathcal{L}}_z$ denote the restriction of $\hat{\mathcal{L}}$ to the fiber 
$U_z:= \hat{\operatorname{pr}}^{-1}(z)$ over $z$, 
and we define $\mathcal{M}'_z 
\subset \Bbb P^*(({\hat{\operatorname{pr}}}_*\hat{\mathcal{L}})_z)$ by
$$
\mathcal{M}'_z\; :=\; \{z\}\times \Phi_{|\hat{\mathcal{L}}_z|}(U_z),
$$
where $\Phi_{|\hat{\mathcal{L}}_z|}$ is the projective embedding 
of $U_z$ associated to the complete linear system $|\hat{\mathcal{L}}_z|$ 
for the line bundle $\hat{\mathcal{L}}_z$.
Then $U$ is viewed as the subset 
$$
U':=\bigcup_{0\neq z\in \Bbb A^1}\mathcal{M}'_z
$$ 
of $\Bbb P^*({\hat{\operatorname{pr}}}_*\hat{\mathcal{L}})$.
Let $\mathcal{M}'$ be the $\Sigma'$-invariant 
subvariety $\bar{U}'$ of $\Bbb P^*({\hat{\operatorname{pr}}}_*\hat{\mathcal{L}})$ obtained 
as the closure of 
$U'$ in $\Bbb P^*({\hat{\operatorname{pr}}}_*\hat{\mathcal{L}})$,
where $U'$ is $\Sigma'$-equivariantly identified with $U$ above. 
The projection
$$
\pi' : \mathcal{M}' \to \Bbb A^1,
$$  
induced by the natural projection 
$\Bbb P^*(\hat{\operatorname{pr}}_*\hat{\mathcal{L}}) 
\to \Bbb A^1$ is a $\Sigma'$-equivariant projective morphism
with a relatively very ample invertible sheaf 
$\mathcal{L}'$, where $\mathcal{L}'$ denotes the restriction of
$\mathcal{O}_{\Bbb P^*(\hat{\operatorname{pr}}_*\hat{\mathcal{L}})}(1)$ 
to $\mathcal{M}'$. 
In view of the $\Sigma'$-action on $\hat{\mathcal{L}}$, note here that $\Sigma'$ acts on $\mathcal{L}'$ 
covering the $\Sigma'$-action on $\mathcal{M}'$.
By replacing $\hat{U}$ by its $\Sigma'$-euivariant modification if necessary, 
we may assume that a $\Sigma'$-equivariant birational surjective morphism
$$
\iota': \hat{U} \to \mathcal{M}'
$$
exists.
We now observe that $\mathcal{M}'\setminus \mathcal{M}'_0 \, =\,  U' 
\,\cong\, U \, =\, \mathcal{M}\setminus \mathcal{M}_0$,
where $\mathcal{M}'_0$ denotes the scheme-theoretic fiber of $\pi'$ over the origin.
Then
$$
(\mathcal{M}'_z, \mathcal{L}'_z)\;\cong\;
(M, L),
\qquad z \in \Bbb A^1 \setminus \{0\}.
$$
Hence $(\mathcal{M}',\mathcal{L}')$ with the $\Bbb C^*$-action is a test configuration for 
the polarized algebraic manifold $(M,L)$, 
where $\Bbb C^* \, (= T_0)$  acts on $\mathcal{M}'$ as the second factor of $\Sigma'
= {\mathcal{S}}'\times \Bbb C^*$.
Then by the notation in the introduction,
$$
\mathcal{S}' \, \subset \,
\mathcal{Q}(\mathcal{M}') \cap G
\, =\, H(\mathcal{M}'). 
$$
Since $Y'$ belongs to the Lie algebra of 
$\mathcal{S}'$,
we see that $Y' \in \frak h (\mathcal{M}')$. 
Moreover $Y'$ is Hamiltonian, and hence
$$
Y' \; \in \; \frak S (\mathcal{M}').
$$
Even for $\mathcal{M}'$ in place of $\mathcal{M}$,
we mean by $X_0$ the generator  
as in the introduction for the $T_0$-action.
Then by setting $X' := X_0 + Y'$, 
we have a generalized test configuration 
$(\mathcal{M}', \mathcal{L}', X')$ of $(M,L)$.

\medskip
For the generalized test configuration $(\mathcal{M}', \mathcal{L}', X')$ defined above,
we set $T':= \{\exp (tX'); t \in \Bbb C\}$, and consider the vector bundle
$$
E'_{m} := \pi'_*({\mathcal{L}'})^{m}, 
\qquad m =1,2,\dots, 
$$
over $\Bbb A^1$.  Note that $\Sigma'$ acts on $E'_m$ covering the natural
$\Sigma'$-action on $\Bbb A^1$.
Then as in (1.2), we have a $\Sigma'$-equivariant trivialization 
for $E'_m$,
$$
E'_m \; \cong \; \Bbb A^1 \times (E'_m)_0.
\leqno{(2.5)}
$$
Consider the 
line bundle $\mathcal{C} := (\iota')^*\mathcal{L}' \otimes (\iota^*\mathcal{L})^{-1}$ 
over $\hat{U}$. 
Then $\mathcal{C}$ is expressible as $\mathcal{O}_{\hat{U}}(\mathcal{D})$
for a divisor $\mathcal{D}$ on $\hat{U}$ satisfying
$\operatorname{Supp}(\mathcal{D})\,\subset\,\hat{\operatorname{pr}}^{-1}(0)$ 
set-theoretically. 
Then the direct image sheaf 
$\Gamma :=\hat{\operatorname{pr}}_*\mathcal{C}$ over $\Bbb A^1$ 
is invertible.
Note that the $T_0$-action on $\hat{U}$ naturally lifts 
to a $T_0$-linearization of $\mathcal{C}$.  Take a 
$T_0$-equivariant trivialization 
$$
\Gamma \; \cong \; \Bbb A^1 \times \Gamma_0,
$$
where $\Gamma_0$ is the fiber of $\Gamma$ over the origin.
Fix an element $0\neq \xi\in\Gamma_0$. Then we have a multiplicative 
algebraic character $\chi : T_0 \to \Bbb C^*$ 
such that $g \cdot \xi = \chi ( g )\, \xi\,$ 
for all $g \in T_0$.
Hence the nonvanishing section $\gamma_{\xi}$ of $\Gamma$ associated to 
$\Bbb A^1 \times \{\xi\}$ satisfies
$$
g \cdot \gamma_{\xi} \; =\; \chi ( g ) \gamma_{\xi},
\qquad g \in T_0.
\leqno{(2.6)}
$$
Now by $H^0(\Bbb A^1, \Gamma )\cong H^0(\hat{U},\mathcal{C})$,
regard $\gamma_{\xi}$ as an element of $H^0(\hat{U},\mathcal{C})$,
and we denote by $(\gamma_{\xi})$ the associated divisor 
on $\hat{U}$.
Then $\mathcal{D}$ above can be chosen as $(\gamma_{\xi})$.
In view the isomorphisms 
$$
\begin{cases}
&H^0(\Bbb A^1, E_m) \cong H^0(\mathcal{M},\mathcal{L}^m) \cong H^0(\hat{U},\iota^*\mathcal{L}^m),\\
&H^0(\Bbb A^1,E'_m) \cong H^0(\mathcal{M}',{\mathcal{L}'}^m) \cong H^0(\hat{U},(\iota' )^*{\mathcal{L}'}^m),
\end{cases} 
$$
the linear isomorphism $H^0(\hat{U},\iota^*\mathcal{L}^m)
\cong H^0(\hat{U},(\iota' )^*{\mathcal{L}'}^m)$, 
$e \leftrightarrow \gamma_{\xi}^{\;m}\cdot e$,
induces a vector bundle isomorphism
$$
E_m \; \cong \; E'_m, 
\qquad m=1,2,\dots,
\leqno{(2.7)}
$$
where in terms of this isomorphism, the $T_0$-actions on $E_m$ and $E'_m$ 
differ just by $\chi^m$, i.e.,
$\Psi^{}_{m,X_0;\mathcal{M}',\mathcal{L}'} (g)\,
=\, \chi (g)^m\,\Psi^{}_{m,X_0;\mathcal{M},\mathcal{L}} (g) $ for all $g\, \in T_0$.
Hence we obtain
$$
\Psi^{\operatorname{SL}}_{m,X_0;\mathcal{M},\mathcal{L}}\; 
=\; \Psi^{\operatorname{SL}}_{m,X_0;\mathcal{M}',\mathcal{L}'},
\leqno{(2.8)}
$$
and we often write them as $\Psi^{\operatorname{SL}}_{m,X_0}$ for 
simplicity. 
Then by the definition of $\mu_{m,s}$ with $z = e^s$, it follows that
$$
(\Pi_m \circ\Psi^{\operatorname{SL}}_{m,X_0;\mathcal{M},\mathcal{L}})(z)\; 
=\; \mu_{m,s}\; =\;
(\Pi_m \circ\Psi^{\operatorname{SL}}_{m,X_0;\mathcal{M}',\mathcal{L}'})(z),
$$
for all $s\in \Bbb C$. 
For the time being, put $m =1$. 
Since $\mathcal{M}_1 = M = \mathcal{M}'_1$ in $\Bbb P^*((E_1)_0)$,
we see from (2.2) that 
$\mathcal{M}_z = \{z\}\times M_s = \mathcal{M}'_z$ for all $z \in \Bbb C^*$,
where $M_s$ is the subvariety $\mu_{1,s}(M)$
in $\Bbb P^*((E_1)_0) \, (= \Bbb P^*((E'_1)_0))$.
Thus
$$
\mathcal{M}\; =\; \mathcal{M}'
\leqno{(2.9)}
$$
as the closure of  
$\mathcal{M}\setminus \mathcal{M}_0 =  \mathcal{M}'\setminus \mathcal{M}'_0 $
in the complex variety $\Bbb A^1\times \Bbb P^*((E_1)^{}_0) = \Bbb A^1\times \Bbb P^*((E'_1)^{}_0)$.

\medskip
Note that both the $T'$-action and the $\mathcal{S}'$-action and on $E'_m$ preserves  
$(E'_m)_0$. We write the associated 
representation on $(E_m)_0$ as
$$
\begin{cases}
&\Psi^{\operatorname{SL}}_{m, X'}
 (=\Psi^{\operatorname{SL}}_{m, X';\mathcal{M}', \mathcal{L}'}) : 
T' \to \operatorname{SL}((E_m)_0)\,
 ( = \operatorname{SL}((E'_m)_0)), \\
&\Psi^{\operatorname{SL}}_{m,\mathcal{S}'} 
(=\Psi^{\operatorname{SL}}_{m, \mathcal{S}';\mathcal{M}', \mathcal{L}'}) \, 
: \,
\mathcal{S}' \to \operatorname{SL}((E_m)_0)\,
 ( = \operatorname{SL}((E'_m)_0)).
\end{cases}\leqno{(2.10)}
$$
By $T_0 \cong \{1\}\times T_0$ and 
$\mathcal{S}' \cong \mathcal{S}' \times \{1\}$, we view
$T_0$, $T'$, $\mathcal{S}'$ as 
subgroups of 
$\Sigma' = \mathcal{S}'\times T_0$. 
Hence $\Psi^{\operatorname{SL}}_{m,X_0}$, $\Psi^{\operatorname{SL}}_{m, X'}$, 
$\Psi^{\operatorname{SL}}_{m,\mathcal{S}'}$
are the restrictions of  
$$
\Psi^{\operatorname{SL}}_{m, \Sigma'}
 (=\Psi^{\operatorname{SL}}_{m, \Sigma';\mathcal{M}', \mathcal{L}'}) : 
 \Sigma' \to \operatorname{SL}((E_m)_0)\,
 ( = \operatorname{SL}((E'_m)_0)).
$$
to these subgroups.
Since the $T_0$-action and the $\mathcal{S}'$-action on $(E_m)_0$ 
commute, 
and since every $\exp (sX')$, $s\in \Bbb C$, in $T'$ is 
expressible as $\exp (sX_0)\cdot\exp (sY')$, 
we obtain
$$
\psi_s\cdot\varphi_s \; =\; 
 \Psi^{\operatorname{SL}}_{m, X'}(\exp (sX')) \; =\;  
 \varphi_s \cdot\psi_s,
 \leqno{(2.11)}
$$
by setting $\varphi_s := \Psi^{\operatorname{SL}}_{m,X_0}(\exp (sX_0))$ 
and $\psi_s := \Psi^{\operatorname{SL}}_{m,\mathcal{S}'}(\exp (sY'))$. 
Note that $\mu_{m, s} = \Pi_m (\varphi_s )$.
For every  $Y \in \frak{S}(\mathcal{M}')\, (= \frak{S}(\mathcal{M}))$,
by setting $X := X_0 + Y$, 
we define
the holomorphic vector field $\mathcal{V}^m = 
\mathcal{V}^{m}(\mathcal{M}', \mathcal{L}', X)$ 
(cf. \cite{M4}) on 
the projective space $\Bbb P^*((E_m)_0)$ by 
$$
\mathcal{V}^m(\mathcal{M}', \mathcal{L}', X)
\; =\;
 (\Pi_m \circ \Psi^{\operatorname{SL}}_{m,\Sigma'})_*(X).
$$
First, we consider the case $Y = Y'$. Then by $X = X'$,
$$
\mathcal{V}^m(\mathcal{M}', \mathcal{L}', X')=
 (\Pi_m \circ \Psi^{\operatorname{SL}}_{m,X'})_*(X').
 $$
 Next let $Y = 0$, so that 
$X = X_0$. Then by (2.8),
$$
\mathcal{V}^m(\mathcal{M}', \mathcal{L}', X_0)
\,  =\, (\Pi_m \circ \Psi^{\operatorname{SL}}_{m,X_0})_*(X_0)\,
=\, \mathcal{V}^m(\mathcal{M}, \mathcal{L}, X_0).
\leqno{(2.12)}
$$

\section{Proof of Main Theorem}

In this section, fixing a polarized algebraic manifold $(M,L)$ 
as in the introduction,
we assume that $c_1(L)$ admits 
a K\"ahler metric 
$$
\omega_{\infty} \; =\; c_1(L, h_{\infty})
$$ 
of constant scalar 
curvature. Here $h_{\infty}$ is a suitably chosen Hermitian metric for $L$.
We here retain the notation in the preceding sections.
Let $\omega_{\operatorname{FS}}$ denote the Fubini-Study metric for 
$\Bbb P^*((E_m)_0)\, ( =\Bbb P^{N_m}(\Bbb C ))$, 
where we identify $\Bbb P^*((E_m)_0 )$ with $\Bbb P^{N_m} (\Bbb C)$ by  
an orthonormal basis for the Hermitian metric $h_{m;0}$ (cf. (1.2)) on $(E_m)_0$.
For each $s \in \Bbb C$, we consider the orthogonal complement $TM_s^{\perp}$ of 
$TM_s$ in 
$T\Bbb P^*((E_m)_0)_{|M_s}$
 in terms of the Fubini-Study metric $\omega_{\operatorname{FS}}$ 
 on $\Bbb P^*((E_m)_0)$.
By writing $h_{m;0}$ simply as $h_m$, and following \cite{M4},
we may choose 
a sequence $h_m$, $m=1,2,\dots$, 
converging in $C^{\infty}$ to $h_{\infty}$ on $M \hookrightarrow \Bbb P^*((E_m)_0)$ 
(cf. (2.1)).
By fixing a Hamiltonian element $W$ in $\frak g$, we consider the associated generalized 
test configuration $(\mathcal{M}', \mathcal{L}', X')$ as in the preceding section,
where $W$ will be specified later.
 Recall that $\mathcal{M}'_z = \{z\}\times M_s$, where $z = e^s$. 
Restricting $\mathcal{V}^{m}(\mathcal{M}', \mathcal{L}', X')$ to $M_s$,
$s \in \Bbb C$,
we can write
$$
\mathcal{V}^m(\mathcal{M}', \mathcal{L}', X')_{|M_s}\; =\;
\mathcal{V}_{TM_{s}}^{m}(\mathcal{M}', \mathcal{L}', X')
+\mathcal{V}_{TM_{s}^{\perp}}^{m}(\mathcal{M}', \mathcal{L}', X'),
\leqno{(3.1)}
$$
where 
$\mathcal{V}_{TM_{s}}^{m}(\mathcal{M}', \mathcal{L}', X')$ 
and $\mathcal{V}_{TM_{s}^{\perp}}^{m}(\mathcal{M}', \mathcal{L}', X')$ 
are smooth sections of $TM_{s}$ and $TM_s^{\perp}$, respectively.
Put $\omega_{s}:= (\mu_{m,s}^*\omega_{\operatorname{FS}}/m)_{|M}$.
For the subspace $\frak r := H^0(M,\mathcal{O}_M(TM))$ of 
$\frak w := H^0(M, C^{\infty}(TM))$, let $\frak r^{\perp}_s$ denote its orthogonal complement 
in terms of the Hermitian pairing
$$
< W_1, W_2>^{}_s\;  := \; 
\int_M \; (W_1, W_2)^{}_{\omega_{s}}\,
\omega_{s}^{\,n},
\qquad W_1, W_2 \in \frak w,
$$
where $(W_1, W_2)^{}_{\omega_{s}}$ is the pointwise Hermitian pairing of $W_1$ and $W_2$ by the K\"ahler metric $\omega_{s}$ on $M$.
Now, as in \cite{M4}, 
$\mathcal{V}_{TM_{s}}^{m} = 
\mathcal{V}_{TM_{s}}^{m}(\mathcal{M}', \mathcal{L}', X')$ is written as 
a sum
$$
\mathcal{V}_{TM_{s}}^{m} \; =\; 
\mathcal{V}_{m,s}^{\circ}+ 
\mathcal{V}_{m,s}^{\bullet},
\leqno{(3.2)}
$$
where $\mathcal{V}_{m,s}^{\circ} = 
\mathcal{V}_{m,s}^{\circ}(\mathcal{M}', \mathcal{L}', X')$ 
and $\mathcal{V}_{m,s}^{\bullet} = 
\mathcal{V}_{m,s}^{\bullet}(\mathcal{M}', \mathcal{L}', X')$ belongs to 
$(\mu_{m,s})_*\frak r$ and $(\mu_{m,s})_*\frak r_s^{\perp}$, respectively.
For $X_0$, $Y'$ and $X'$ regarded as elements of 
the Lie algebra of $\Sigma'$, we consider the $\Sigma'$-action 
on $(E'_m)_0\, (= (E_m)_0)$. Then by $X_0^{(m)}$, ${Y'}^{(m)}$ and ${X'}^{(m)}$,
we mean the holomorphic vector fields on the projective space 
$\Bbb P^*((E'_m)_0)\, (= \Bbb P^*((E_m)_0))$
induced by the infinitesimal actions of $X_0$, 
$Y'$ and $X'$, respectively.
We now observe that  ${{X'}^{(m)}} = \mathcal{V}^{m}(\mathcal{M}', \mathcal{L}', X')$ 
and that ${X_0^{(m)}} = \mathcal{V}^{m}(\mathcal{M}', \mathcal{L}', X_0)$.
Hence by ${{X'}^{(m)}} = {X_0^{(m)}} +{ {Y'}^{(m)}}$ 
and (2.12), it follows that
$$
\begin{cases}
\;\;&\mathcal{V}^{m}(\mathcal{M}', \mathcal{L}', X') - { {Y'}^{(m)}}
=\; {X_0^{(m)}}\;  \\
&=\; \mathcal{V}^{m}(\mathcal{M}', \mathcal{L}', X_0)\;
=\; \mathcal{V}^{m}(\mathcal{M}, \mathcal{L}, X_0).
\end{cases}
\leqno{(3.3)}
$$
By (2.1), we identify $M$ with $\mathcal{M}'_1\,( = \mathcal{M}_1)$ 
sitting in $\Bbb P^*((E'_m)_0)\, (= \Bbb P^*((E_m)_0))$. 
Let us now put
$$
A_m\,:=\, {Y'^{(m)}}_{|M} \quad \text{and} \quad
B_m\, :=\, X_0^{(m)}{}_{|M}.
$$
Then $A_m$ coincides with the holomorphic vector field $Y' = -W$ on $M$, 
and by (3.3), $B_m$ coincides with 
$\mathcal{V}^m(\mathcal{M}, \mathcal{L}, X_0 )_{|M}$.
The element $\varphi_s$ in $\operatorname{SL}((E_m)_0)\, (=\operatorname{SL}((E'_m)_0))$ (cf.~(2.11))
induces a projective linear transformation on $\Bbb P^*((E'_m)_0)\, (= \Bbb P^*((E_m)_0))$,
taking $M$ to $M_s$. This mapping of $M$ onto $M_s$ 
is written also as $\varphi_s$ by abuse of terminology.
Then for each $s \in \Bbb C$, we can write
\begin{align*}
{ {Y'}^{(m)}}_{|M_s}  \; &=\; (\varphi_s)_* A_m\; =\; -(\varphi_s)_* W,
\tag{3.4}
\\
X_0^{(m)}{}_{|M_s} 
 \; &=\; 
 (\varphi_s)_*B_m \; =\; 
 (\varphi_s)_* \{\mathcal{V}_m(\mathcal{M}, \mathcal{L}, X_0 )_{|M}\},
\tag{3.5}
\end{align*} 
where (3.4) follows from the commutativity of the $\mathcal{S}'$-action and 
$T_0$-action. Now by (3.4), 
${ {Y'}^{(m)}}_{|M_s} $ is a holomorphic vector 
field on $M_s$, 
and hence by (3.3), we obtain
$$
\begin{cases}
&\mathcal{V}^{m}_{TM^{\perp}_s}(\mathcal{M}', \mathcal{L}', X') \,=\, 
\mathcal{V}^{m}_{TM^{\perp}_s}(\mathcal{M}, \mathcal{L}, X_0),\\
&\mathcal{V}^{\bullet}_{m,s}(\mathcal{M}', \mathcal{L}',X')\;\;\, =\;\;\, 
\mathcal{V}^{\bullet}_{m,s}(\mathcal{M}, \mathcal{L}, X_0),
\end{cases}\leqno{(3.6)}
$$
where we used the identification
$\mathcal{M}'_{z}= \{z\}\times M_s =\mathcal{M}_{z}$ 
(cf. (2.9)) with $z = e^s$.
Then by (3.3), (3.4) and (3.5),
${\mathcal{V}^m(\mathcal{M}', \mathcal{L}', X')}_{|M_s}\, (= {X'^{(m)}}_{|M_s})$ 
 is expressible as 
$$
{\mathcal{V}^m(\mathcal{M}', \mathcal{L}', X')}_{|M_s}\;  =\; 
(\varphi_s)_*\{-W + \mathcal{V}^m(\mathcal{M}, \mathcal{L}, X_0 )_{|M}\}.
 \leqno{(3.7)}
$$
Then for a subsequence $\{m(j)\,;\,j =1,2,\dots\}$ of $\{m =1,2,\dots \}$, 
and a sequence $\{s(j)\in \Bbb R\,;\,j=1,2,\dots\}$  
as in Section 5 in \cite{M4} such that
$$
- \,C_0\,(\log m(j) )\, q (j) \leq s(j) \leq 0
\quad\text{ and } \quad q (j) = 1/m (j),
$$
we obtain (cf. (5.19) and (5.20) of \cite{M4})
\begin{align*}
\int_M |\mathcal{V}_{TM^{\perp}}(j)|^2_{\omega (j)}\, \omega (j)^n &= O 
\left (\frac{q(j)}{\log m(j)}\right ),
\tag{3.8}\\
\int_M |\mathcal{V}^{\bullet}(j)|^2_{\omega (j)}\, \omega (j)^n \;\;&= O 
\left (\frac{1}{\log m(j)}\right ),
\tag{3.9}
\end{align*}
where $C_0$ is a sufficiently small positive real constant independent of $j$, and
we define $\omega (j)$, $N(j)$, $\mathcal{V}_{TM^{\perp}}(j)$, 
$\mathcal{V}^{\bullet}(j)$, $\mathcal{V}^{\circ}(j)$ by
$$
\begin{cases}
\;\;\omega (j) &:= \,\{ q(j) \,\mu_j^*\omega_{\operatorname{FS}}\}_{|M},\;\;
N(j):= N_{m(j)},\\
\;\;\mathcal{V}_{TM^{\perp}}(j) &:=
(\mu_j^{-1})_*\mathcal{V}^{m(j)}_{TM^{\perp}_{s(j)}}
(\mathcal{M}, \mathcal{L}, X_0)
=(\mu_j^{-1})_*\mathcal{V}^{m(j)}_{TM^{\perp}_{s(j)}}
(\mathcal{M}', \mathcal{L}', X') ,\\
\;\;\mathcal{V}^{\bullet}(j)& := (\mu_j^{-1})_*\mathcal{V}^{\bullet}_{m(j),s(j)}(\mathcal{M}, \mathcal{L}, X_0)
= (\mu_j^{-1})_*\mathcal{V}^{\bullet}_{m(j),s(j)}(\mathcal{M}', \mathcal{L}', X'),\\
\;\;\mathcal{V}^{\circ}(j)& := (\mu_j^{-1})_*\mathcal{V}^{\circ}_{m(j),s(j)}(\mathcal{M}, \mathcal{L}, X_0),
\end{cases}
$$
with $\mu_j := \mu_{m(j), s(j)}$. Put $\mathcal{V}^{\circ}(j)' := (\mu_j^{-1})_*\mathcal{V}^{\circ}_{m(j),s(j)}(\mathcal{M}', \mathcal{L}', X')$.  
In view of $\mu_{m,s} = \Pi_m (\varphi_{s})$, we can write $(\varphi_s)_*$ in (3.7) 
as $(\mu_{m,s})_*$.
Then by (3.7) applied to $s(j)$, it follows from (3.6) that
$$
\mathcal{V}^{\circ}(j)' \; =\; \mathcal{V}^{\circ}(j) - W.
\leqno{(3.10)}
$$
In (5.1) of \cite{M4}, the basis $\{\tau_1, \tau_2, \dots, \tau_{N(j)}\}$ 
for $(E_{m(j)})_0$
can be chosen in such a way that, in terms of this basis, the $\Sigma'$-action 
(as well as the $T_0$-action)
on $(E_m)_0$ is also diagonalizable.
Then $\eta (j) = \eta (j) (\mathcal{M}, \mathcal{L}, X_0) \in C^{\infty}(M)_{\Bbb R}$ 
as in \cite{M4}
satisfies
$$
i_{\mathcal{V}(j)}\omega (j)\;  = \; \bar{\partial}\eta (j), 
$$
where $\mathcal{V}(j):= 
\{(\mu_j^{-1})_*\mathcal{V}^{m(j)}(\mathcal{M}, \mathcal{L}, X_0 )\}_{|M}\,( = \mathcal{V}^{m(j)}(\mathcal{M}, \mathcal{L}, X_0 )_{|M})$.
By (3.7), we see that 
$\eta' (j) := \eta (j) - f_{\omega (j), W} \in C^{\infty}(M)_{\Bbb R}$  (cf. (1.1))
satisfies
$$
i_{\mathcal{V}'(j)}\omega (j) \; =\; \bar{\partial} \eta' (j),
\leqno{(3.11)}
$$
where $\mathcal{V}'(j):
=\{(\mu_j^{-1})^*\mathcal{V}^m(\mathcal{M}', \mathcal{L}', X')\}_{|M}$. 
Recall that, as $j \to \infty$, 
the K\"ahler metric $\omega (j )$ converges in $C^{\infty}$ to the metric $\omega_{\infty}$ 
in $c_1(L)_{\Bbb R}$ of constant scalar curvature.  Moreover,
for some $\eta_{\infty}\in C^{\infty}(M)_{\Bbb R}$,
$$
\eta (j) \to \eta_{\infty} \text{ in $L^2(M, \omega_{\infty}^n)$}, 
\qquad \text{ as $j\to \infty$},
\leqno{(3.12)}
$$
where $\eta_{\infty}$ is a `Hamiltonian function' for a holomorphic vector field 
on $M$ (cf. \cite{M4}). 
Then $W\in \frak g$ is specified as follows.
Following \cite{M4}, by using the notation in (1.1), let
$W$ be the unique element of $\frak g$ such that  
$$
f_{\omega_{\infty}, W} = \eta_{\infty}.
\leqno{(3.13)}
$$
Since $\omega (j )$ converges to $\omega_{\infty}$ in $C^{\infty}$,
by (3.12) and (3.13), 
the definition of $\eta'(j)$ implies the convergence
$$
\eta'(j) \to 0 \text{ in $L^2(M)$}, \qquad \text{ as $j\to \infty$}.
\leqno{(3.14)}
$$
By the identification of $(E_{m(j)})_1 \,(= (E'_{m(j)})_1) $ with 
$(E_{m(j)})_0\, (= (E_{m(j)})_0)$ (cf. (2.5)),
we consider the basis $\{\tau'_1, \tau'_2, \dots, \tau'_{N(j)}\}$ for
$$
 H^0(M,\mathcal{O}_M(L^{m(j)})) = (E_{m(j)})_1
$$ 
corresponding to the basis $\{\tau_1, \tau_2, \dots, \tau_{N(j)}\}$ for 
$(E_{m(j)})_0$.
By the representation in the first line of (2.10),
we can write
$$
X'\cdot \tau_{\alpha} \; =\; e'_{\alpha}(j)\, \tau_{\alpha}, 
\qquad \alpha =1,2,\dots, N(j),
\leqno{(3.15)}
$$
for some real numbers $e'_{\alpha}(j)$ with $\Sigma_{\alpha =1}^{N(j)}\, 
e'_{\alpha}(j) = 0$.
By setting
$$
\xi (j) 
:= \frac{\Sigma_{\alpha =1}^{N(j)}e'_{\alpha}(j)|z_{\alpha}|^2}{ m(j)\,
\Sigma_{\alpha =1}^{N(j)}|z_{\alpha}|^2}
\quad \text{and}\quad
\bar{\xi} (j) := \; \frac{\int_{M_{s(j)}} \xi (j)\, \{ q (j)\, \omega_{\operatorname{FS}} \}^n}
{\int_{M_{s(j)}} \{q (j)\, \omega_{\operatorname{FS}} \}^n},
$$
we consider $e_{\alpha}(j) := e'_{\alpha}(j) - m(j)\bar{\xi} (j)$, where $\Bbb P^*((E_{m(j)})_0)$ 
is identified with the projective space 
$\Bbb P^{N(j)}(\Bbb C ) = \{ (z_1:z_2:\dots:z_{N(j)})\}$ 
by the basis $\{\tau_1, \tau_2, \dots, \tau_{N(j)}\}$.  
Hereafter, let $C_k$, $k$=1,2,\dots, denote positive 
real constants independent of $j$ and $\alpha$.
Since $|e_{\alpha}'(j)|\leq C_1 m(j)$ for all $\alpha$ and $j$,
the inequality $|\xi (j) | \leq C_1$ holds for all $j$, 
so that  
$$
|\bar{\xi} (j) | \leq C_1 \quad\text{ for all $j$}.
\leqno{(3.16)}
$$
In particular $|e_{\alpha}(j)|\leq C_2\, m(j)$ for all $j$.
Put $h(j):= h_{m(j)}$. In view of the definition of $\eta (j)$ in \cite{M4}, 
we can write $\eta' (j)$ in (3.11) as
$$
\eta'(j)\; = \;  
\frac{\Sigma_{\alpha = 1}^{N(j)}\, e_{\alpha}(j)\, |\tau'_{\alpha}|^2_{h(j)}
\exp \{2s(j) e_{\alpha}(j)\}}
{m(j) \Sigma_{\alpha = 1}^{N(j)}\, |\tau'_{\alpha}|^2_{h(j)}
\exp \{2s(j) e_{\alpha}(j)\} },
$$
since $\int_M \eta' (j) \omega (j )^n =0$.
Note that the functions $\eta' (j)$, $j=1,2,\dots$, are uniformly bounded on $M$.
We further define uniformly bounded functions 
$\zeta'(j)$, $j =1,2,\dots$, on $M$ by
$$
\zeta' (j) \; :=\; 
\frac{\Sigma_{\alpha = 1}^{N(j)}\, e_{\alpha}(j)^{\,2}\, |\tau'_{\alpha}|^2_{h(j)}
\exp \{2s(j) e_{\alpha}(j)\}}
{m(j)^2 \Sigma_{\alpha = 1}^{N(j)}\, |\tau'_{\alpha}|^2_{h(j)}
\exp \{2s(j) e_{\alpha}(j)\} }.
$$
Define bounded sequences $\{\lambda_j\}$, $\{\kappa_j\}$ of real numbers by 
setting
$$
\lambda_j := \int_M \eta' (j)^2\, \omega (j)^n\,\geq \, 0\,
\quad\text{ and }\quad \kappa_j := \int_M \zeta' (j) \, \omega (j)^n\,\geq \,0.
$$
If necessary, we replace $\{\lambda_j\}$,
$\{\kappa_j\}$  by respective subsequences with a common sequence of indices.
Then we may assume that 
both $\{\lambda_j\}$ and $\{\kappa_j\}$ converge.
Put 
$$
\lambda_{\infty} := \lim_{j\to\infty} \lambda_j
\quad\text{ and }\quad
\kappa_{\infty} := \lim_{j \to \infty} \kappa_j.
$$
By the Cauchy-Schwarz inequality, 
we have $\eta'(j)^2 \leq \zeta'(j)$ and hence
$\kappa_j \geq \lambda_j$ for all $j$. 
In particular $\kappa_{\infty} \geq \lambda_{\infty}$. 
We now claim that there exists a positive constant $C$ independent of $j$ such that
$$
\text{$0 \leq \kappa_j - \lambda_j \leq C/m(j)\;$ and hence $\;\kappa_{\infty} =
 \lambda_{\infty}$}.
\leqno{(3.17)}
$$
{\em Proof.} 
For each $j =1,2,\dots$, we define $I_j^{\circ}$ and ${I_j^{\circ}}'$ by
$$
\begin{cases}
\; I_j^{\circ} &:= \int_M |\mathcal{V}^{\circ} (j)|^2_{\omega (j)} \, \omega (j)^n, \\
\; {I_j^{\circ}}'&:= \int_M |\mathcal{V}^{\circ} (j)'|^2_{\omega (j)} \, \omega (j)^n.
\end{cases}
$$
Since it suffices to consider Case 1 of Section 5 in \cite{M4},
we may assume that $I_j^{\circ}$,
$j$=1,2,\dots, are bounded. Then by (3.10), 
${I_j^{\circ}}'$,
$j$=1,2,\dots, are also bounded. Since
$I'_j\, :=\,\int_M |\mathcal{V}'(j)|^2_{\omega (j)}\,\omega (j)^n$ is written as 
$$
I'_j \; =\; {I_j^{\circ}}' + \int_M |\mathcal{V}^{\bullet}(j)|^2_{\omega (j)}\, \omega (j)^n
+ \int_M |\mathcal{V}_{TM^{\perp}}(j)|^2_{\omega (j)}\, \omega (j)^n, 
$$
we see from (3.8) and (3.9) that $I_j'$, $j$=1,2,\dots, are also bounded,
 i.e., $0\leq I_j \leq C$ for some $C$ as above.
On the other hand, it is easily checked that $I'_j = m(j)(\kappa_j - \lambda_j)$.  
Hence $0 \leq \kappa_j - \lambda_j \leq C/m(j)$. By letting $j\to \infty$,
we obtain
 $\kappa_{\infty} = \lambda_{\infty}$, as required.
 \qed

\medskip
Now by (3.14), $\lambda_{\infty} = 0$, and hence from (3.17), 
it follows that 
$$
\kappa_{\infty} = 0.
\leqno{(3.18)}
$$
If $\Psi_{m,X'}^{\operatorname{SL}}$ in (2.10)
is nontrivial for $m =1$,  then from the next section, we see that 
$\kappa_{\infty} > 0$  in contradiction to (3.18).
Hence, $\Psi_{1,X'}^{\operatorname{SL}}$
is trivial. 
Then the argument leading to (2.9) shows that
$\mathcal{M}'_z = \{z\}\times M_s$, where $M_s$ is the subvariety
$\mu_{1,s} (M) = (\varphi_s \psi_s) (M)$ 
in $\Bbb P^*((E_1)_0)$.
Thus by the triviality of $\Psi_{1,X'}^{\operatorname{SL}}$, 
we see from (2.11) that the generalized test configuration 
$(\mathcal{M}', \mathcal{L}', X')$  is trivial, i.e., $\mathcal{M}' = \Bbb A^1 \times M$ 
with $T'$ acting on the second factor trivially. Therefore by (2.8) and (2.11),
we conclude that the original test configuration $(\mathcal{M},\mathcal{L})$ 
is a product configuration, i.e., $\mathcal{M} = \Bbb A^1 \times M$, 
as required.
\qed

\section{$\kappa_{\infty} > 0$ if
$\Psi_{1,X'}^{\operatorname{SL}}$ is nontrivial}

In this section, we consider the set  $\Delta^+_j$
of all $\alpha\in \{1,2,\dots, N(j)\}$  such that $e_{\alpha}(j) \geq 0$.
Similarly, let $\Delta^-_j$ be the set of all 
$\alpha\in \{1,2,\dots, N(j)\}$ such that $e_{\alpha}(j) <  0$.
Put
\begin{align*}
&\nu_{\alpha}(j)\; : = \int_M  \frac{ |\tau'_{\alpha}|^2_{h(j)}
\exp \{2s(j) e_{\alpha}(j) \}}{\Sigma_{\alpha = 1}^{N(j)}\, |\tau'_{\alpha}|^2_{h(j)}
\exp \{2s(j) e_{\alpha}(j)\}} \, \omega (j)^n,\\
&\hat{e}_{\alpha}(j) \; := \; e_{\alpha}(j)/m(j),
\end{align*}
for $\alpha = 1,2,\dots, N(j)$. 
Then $\Sigma_{\alpha =1}^{N(j)}\,\nu_{\alpha}(j) = 1$
and $|\hat{e}_{\alpha}(j) | \leq C_2$.
Since $\int_M \eta'(j)\omega (j)^n =0$ and 
$\kappa_j = \Sigma_{\alpha = 1}^{N(j)}\, \nu_{\alpha}(j)  e_{\alpha}(j)^2 /m(j)^2$, 
we obtain
\begin{align*}
\Sigma_{\alpha \in \Delta^+_j}\; \nu_{\alpha}(j)\, \hat{e}_{\alpha}(j) \; &=\;
-\Sigma_{\alpha \in \Delta^-_j}\; \nu_{\alpha}(j)\, \hat{e}_{\alpha}(j),
\tag{4.1}\\
\kappa_j \; &=\; \Sigma_{\alpha =1}^{N(j)}
\,\nu_{\alpha}(j)  \,\hat{e}_{\alpha}(j)^2.
\tag{4.2}
\end{align*}
Let $r_j$ be the left-hand side (= right-hand side) of (4.1). 
From now on until the end of this section, we assume that 
$\Psi_{1,X'}^{\operatorname{SL}}$ is nontrivial.
Now we claim the following inequality: 
$$
r_{\infty}\, :=\, \varlimsup_{j \to \infty} \;r_j \; >\;  0.
\leqno{(4.3)}
$$
On the other hand, by $0\leq \Sigma_{\alpha\in \Delta_j^+}\,\nu_{\alpha}(j)\leq 
1$,
the Cauchy-Schwarz inequality together with (4.2) implies that
\begin{align*}
r_j^{\,2}\, &=\, 
\{\Sigma_{\alpha \in \Delta^+_j}\; \nu_{\alpha}(j)\, \hat{e}_{\alpha}(j)\}^2
 \leq 
\{\Sigma_{\alpha \in \Delta_j^+}\,\nu_{\alpha}(j)\}
\{\Sigma_{\alpha \in \Delta_j^+}\,\nu_{\alpha}(j)\,\hat{e}_{\alpha}(j)^2\}\\
&\leq \,\Sigma_{\alpha \in \Delta_j^+}\,\nu_{\alpha}(j)\,\hat{e}_{\alpha}(j)^2
\leq \kappa_j.
\end{align*}
We here let $j\to\infty$. Hence, once (4.3) is proved, it follows that
$$
0 \; <\; r_{\infty}^{\,2}\; \leq \; \kappa_{\infty},
$$
as required.
Thus the proof of $\kappa_{\infty} > 0$ is reduced to showing (4.3). 

\medskip\noindent
{\em Proof of $(4.3)$.}  For each $\, \theta\in \Bbb R$ with $0\leq \theta \leq C_2$, 
we consider the set $\Delta_j^{\theta,+}$ of all 
$\alpha\in \{1,2,\dots,N(j)\}$ such that 
$\hat{e}_{\alpha}(j) \geq \theta$. Consider also the set $\Delta_j^{\theta,-}$
 of all 
$\alpha\in \{1,2,\dots,N(j)\}$ such that $\hat{e}_{\alpha}(j) < \theta$.
For each 
$\alpha \in \{1,2,\dots,N(j)\}$, we put 
$\hat{\epsilon}_{\alpha}(j):= \hat{e}_{\alpha}(j) - \theta$ 
and $\epsilon_{\alpha}(j):= m(j) \,\hat{\epsilon}_{\alpha}(j) = e_{\alpha}(j)\, -\, m(j)\,\theta$.
We now define nonnegative real-valued functions 
$V^{\theta,+}(j)$, $V^{\theta,-}(j)$, 
$Y^{\theta,+}$, $Y^{\theta,-}$, $Z^{\theta,+}(j)$, $Z^{\theta,-}(j)$ on $M$ by
\begin{align*}
V^{\theta,+}(j) &:= \; \Sigma^{}_{\alpha\in \Delta^{\theta,+}_j}\,
\hat{\epsilon}_{\alpha}(j)^2\, |\tau'_{\alpha}|^2_{h(j)}
\exp \{2s(j)\, \epsilon_{\alpha}(j) \},\\
V^{\theta,-}(j) &:= \; \Sigma^{}_{\alpha\in \Delta^{\theta,-}_j}\,
\hat{\epsilon}_{\alpha}(j)^2\, |\tau'_{\alpha}|^2_{h(j)}
\exp \{2s(j)\, \epsilon_{\alpha}(j) \},\\
Y^{\theta,+}(j) &:=\; \Sigma^{}_{\alpha\in \Delta^{\theta,+}_j}\,
\hat{\epsilon}_{\alpha}(j)\,|\tau'_{\alpha}|^2_{h(j)}
\exp \{2s(j)\, \epsilon_{\alpha}(j) \},\\
Y^{\theta,-}(j) &:=\; \Sigma^{}_{\alpha\in \Delta^{\theta,-}_j}\,
\{-\hat{\epsilon}_{\alpha}(j)\}\,|\tau'_{\alpha}|^2_{h(j)}
\exp \{2s(j)\, \epsilon_{\alpha}(j) \}.\\
Z^{\theta,+}(j) &:=\; \Sigma^{}_{\alpha\in \Delta^{\theta,+}_j}\,
 |\tau'_{\alpha}|^2_{h(j)}
\exp \{2s(j)\, \epsilon_{\alpha}(j)\},\\
Z^{\theta,-}(j) &:=\; \Sigma^{}_{\alpha\in \Delta^{\theta,-}_j}\,
 |\tau'_{\alpha}|^2_{h(j)}
\exp \{2s(j)\, \epsilon_{\alpha}(j)\},\\
Z(j) \,\;\;\; &:=\; \Sigma_{\alpha = 1}^{N(j)}\, |\tau'_{\alpha}|^2_{h(j)}
\exp \{2s(j)\, \epsilon_{\alpha}(j)\}\; =\; Z^{\theta,+}(j) + Z^{\theta,-}(j).
\end{align*}
Then by the Cauchy-Schwarz inequality, we have
\begin{align*}
\Theta_+& := \frac{V^{\theta,+}(j)Z^{\theta,+}(j) 
- Y^{\theta,+}(j)^2}{Z(j)^2} \geq 0,\\
\Theta_-& := \frac{V^{\theta,-}(j)Z^{\theta,-}(j) 
- Y^{\theta,-}(j)^2}{Z(j)^2} \geq 0.
\end{align*} 
Hence by
$\kappa_j - \lambda_j = \int_M \{\zeta'(j) - \eta'(j)^2\} \omega (j)^n 
\leq C/m(j)$ (cf. (3.17)), and also by 
\begin{align*}
0 &\leq \zeta'(j) - \eta'(j)^2 \\
&= \Theta_+\, +\, \Theta_-\, + \,
\frac{V^{\theta,+}(j)Z^{\theta,-}(j)+ V^{\theta,-}(j)Z^{\theta,+}(j)
+ 2Y^{\theta,+}(j)Y^{\theta,-}(j)}{Z(j)^2},
\end{align*}
we see the following inequality:
\begin{align*}
&0 \leq \int_M \frac{V^{\theta,+}(j)Z^{\theta,-}(j)}{Z(j)^2} \omega (j)^n \leq \frac{C}{m(j)},
\tag{4.4}\\
&0 \leq \int_M \frac{V^{\theta,-}(j)Z^{\theta,+}(j)}{Z(j)^2} \omega (j)^n \leq \frac{C}{m(j)}.
\tag{4.5}
\end{align*}
We further define nonnegative real-valued functions $R^{\theta,+}(j)$, 
$R^{\theta,-}(j)$, $R(j)$, $S^{\theta,+}(j)$, $S^{\theta,-}(j)$, $U^{\theta,+}(j)$, $U^{\theta,-}(j)$ on $M$ by
\begin{align*}
R^{\theta,+}(j) &:=\; \Sigma^{}_{\alpha\in \Delta^{\theta,+}_j}\,
|\tau'_{\alpha}|^2_{h(j)},\\
R^{\theta,-}(j) &:=\; \Sigma^{}_{\alpha\in \Delta^{\theta,-}_j}\,
|\tau'_{\alpha}|^2_{h(j)},\\
R(j)\;\;\;\; &:=\; \Sigma_{\alpha = 1}^{N(j)}|\, \tau'_{\alpha}|^2_{h(j)}\; 
=\; R^{\theta,+}(j) + R^{\theta,-}(j),\\
S^{\theta,+}(j) &:= \; \Sigma^{}_{\alpha\in \Delta^{\theta,+}_j}\,
\hat{\epsilon}_{\alpha}(j)\, |\tau'_{\alpha}|^2_{h(j)},\\
S^{\theta,-}(j) &:= \; \Sigma^{}_{\alpha\in \Delta^{\theta,-}_j}\,
\{-\hat{\epsilon}_{\alpha}(j)\}\, |\tau'_{\alpha}|^2_{h(j)},\\
U^{\theta,+}(j) &:= \; \Sigma^{}_{\alpha\in \Delta^{\theta,+}_j}\,
\hat{\epsilon}_{\alpha}(j)^2\, |\tau'_{\alpha}|^2_{h(j)},\\
U^{\theta,-}(j) &:= \; \Sigma^{}_{\alpha\in \Delta^{\theta,-}_j}\,
\hat{\epsilon}_{\alpha}(j)^2\, |\tau'_{\alpha}|^2_{h(j)}.
\end{align*}
Now by $-C_0(\log m(j))q(j) \leq s(j) \leq 0$ and $|\hat{\epsilon}_{\alpha} (j) | \leq 2 C_2$,
we have
$$
m(j)^{-4C_0C_2}\leq \exp\{2 s(j)\,\epsilon_{\alpha}(j)\} \; \leq \; m(j)^{4C_0C_2},
$$
Since $C_0$ is sufficiently small, we may assume $C_0C_2 \leq  \varepsilon /16$
for a positive constant $\varepsilon \ll 1$ independent of $j$ and $\theta$. 
Hence
$$
\begin{cases}
&0\leq m(j)^{-\varepsilon /4}Z(j) \leq  R(j)\leq m(j)^{\varepsilon /4}Z(j) ,\\
&0\leq m(j)^{-\varepsilon /4}S^{\theta,+}(j) \leq  Y^{\theta,+}(j)
\leq S^{\theta,+}(j),\\
&0\leq S^{\theta,-}(j) \leq  Y^{\theta,-}(j),\\
&0\leq m(j)^{-\varepsilon /4}U^{\theta,+}(j) \leq  V^{\theta,+}(j),\\
&0\leq U^{\theta,-}(j) \leq  V^{\theta,-}(j),\\
&0\leq Z^{\theta,+}(j) \leq  R^{\theta,+}(j)
\leq m(j)^{\varepsilon /4}Z^{\theta,+}(j),\\
&0\leq m(j)^{-\varepsilon /4}Z^{\theta,-}(j) \leq  R^{\theta,-}(j)\leq Z^{\theta,-}(j).
\end{cases} \leqno{(4.6)}
$$
Then by (4.4) and (4.5),
\begin{align*}
&0\leq \int_M \frac{U^{\theta,+}(j)R^{\theta,-}(j)}{R(j)^2}\omega (j)^n 
\leq \frac{C}{m(j)^{1-\varepsilon}}, \tag{4.7}\\
&0\leq \int_M \frac{U^{\theta,-}(j)R^{\theta,+}(j)}{R(j)^2}\omega (j)^n 
\leq \frac{C}{m(j)^{1-\varepsilon}}, \tag{4.8}
\end{align*}
while by $S^{\theta,+}(j)^2 \leq U^{\theta,+}(j)R^{\theta,+}(j)$ 
and $S^{\theta,-}(j)^2 \leq U^{\theta,-}(j)R^{\theta,-}(j)$, we have the
inequalities
\begin{align*}
&\left (\frac{S^{\theta,+}(j) R^{\theta,-}(j)}{R(j)^2} \right )^2 \leq
\frac{S^{\theta,+}(j)^2 R^{\theta,-}(j)}{R^{\theta,+}(j)R(j)^2} 
\leq \frac{U^{\theta,+}(j)R^{\theta,-}(j)}{R(j)^2},
\tag{4.9} \\
&\left (\frac{S^{\theta,-}(j) R^{\theta,+}(j)}{R(j)^2} \right )^2 \leq
\frac{S^{\theta,-}(j)^2 R^{\theta,+}(j)}{R^{\theta,-}(j)R(j)^2} \leq 
\frac{U^{\theta,-}(j)R^{\theta,+}(j)}{R(j)^2}.
\tag{4.10}
\end{align*}
It now follows from (4.7) and (4.9) that
$$
\int_M \frac{S^{\theta,+}(j) R^{\theta,-}(j)}{R(j)^2} \omega (j)^n \leq 
 \frac{C_3}{\; m(j)^{(1-\varepsilon )/2} },
 \leqno{(4.11)}
$$
where $C_k$, $k=3,4,\dots$, are positive constants independent of $j$, $\alpha$ 
and also $\theta$.
Similarly by (4.8) and (4.10), we obtain
$$
\int_M \frac{S^{\theta,-}(j) R^{\theta,+}(j)}{R(j)^2} \omega (j)^n \leq 
 \frac{C_4}{\; m(j)^{(1-\varepsilon)/2} }.
 \leqno{(4.12)}
$$
Put $\delta:= 1/\log m(j)$.
Consider the set $\Delta_j^{\theta,\delta}$ of all $\alpha \in \{1,2,\dots,N(j)\}$ 
such that $\theta\leq \hat{e}_{\alpha}(j) < \theta + \delta$, 
i.e., $0 \leq \hat{\epsilon}_{\alpha}(j) \leq \delta$.
By setting
$$
\begin{cases}
&Y^{\theta,\delta}:= \Sigma_{\alpha\in \Delta_j^{\theta,\delta}}\;
\hat{e}_{\alpha}(j) |\tau'_{\alpha}|^2_{h(j)}
\exp\{2s(j)\epsilon_{\alpha}(j)\},\\
&\tilde{Y}^{\theta,\delta}:= \Sigma_{\alpha\in \Delta_j^{\theta,\delta}}\;
\hat{\epsilon}_{\alpha}(j) |\tau'_{\alpha}|^2_{h(j)}
\exp\{2s(j)\epsilon_{\alpha}(j)\},
\end{cases}
$$
we consider the nonnegative real-valued functions
$\rho^{\theta,\delta}(j)$, $\tilde{\rho}^{\theta,\delta}(j)$, $\rho^{\theta,-}(j)$, 
$\tilde{\rho}^{\theta,-}(j)$, $j=1,2,\dots$, 
on $M$ defined by
\begin{align*}
&\rho^{\theta,\delta}(j) :=  \frac{Y^{\theta,\delta}(j)}{Z(j)},\;
\tilde{\rho}^{\theta,\delta}(j) :=\;  
\frac{\theta \,\delta^{-1} \tilde{Y}^{\theta,\delta}(j)}{R(j)},\\
&\rho^{\theta,-}(j) :=  \frac{Y^{\theta,-}(j)}{Z(j)}, \;
\tilde{\rho}^{\theta,-}(j) :=  \frac{S^{\theta,-}(j) }{R(j)}. 
\end{align*}
We first study the function $\rho^{\theta,\delta}(j)$ and integrate it over $M$.
Note that $\theta \,\delta^{-1}\hat{\epsilon}_{\alpha}(j) \leq \theta \leq 
\hat{e}_{\alpha}(j)$ for all $\alpha \in \Delta_j^{\theta,\delta}$.
Hence
\begin{align*}
&{\rho}^{\theta,\delta}(j) - \tilde{\rho}^{\theta,\delta}(j) 
\;\geq \;
\theta \,\delta^{-1}
\left ( \frac{\tilde{Y}^{\theta,\delta}(j)}{Z(j)} -  \frac{ \tilde{Y}^{\theta,\delta}(j)}{R(j)}
\right ) \\
&\geq \; 
\theta \,\delta^{-1} \frac{\tilde{Y}^{\theta,\delta}(j)\{R^{\theta,+}(j)-Z^{\theta,+}(j) \} 
+\tilde{Y}^{\theta,\delta}(j)\{ R^{\theta,-}(j) - Z^{\theta,-}(j)\} }{Z(j) R(j)}\\
& \geq \; 
-\,\theta \,\delta^{-1}
\frac{\tilde{Y}^{\theta,\delta}(j)Z^{\theta,-}(j)  }{Z(j) R(j)}\;
\geq \; 
-\,\theta \,\delta^{-1}
\frac{{Y}^{\theta,+}(j)Z^{\theta,-}(j) }{Z(j) R(j)}.
\end{align*}
Integrating this over $M$, by (4.6) and (4.11), we obtain
$$
 \int_M \rho^{\theta, \delta}(j) \omega (j)^n
\geq  \int_M \tilde{\rho}^{\theta,\delta}(j) \omega (j)^n - 
\theta\,\delta^{-1}\frac{C_5}{m(j)^{(1/2)-\varepsilon}}.
$$
Since $\Sigma_{\alpha\in \Delta_j^{\theta, \delta}}\; 
\nu_{\alpha}(j)\,\hat{e}_{\alpha}(j)
= \int_M \rho^{\theta, \delta}(j) \omega (j)^n$, it follows that
$$
\Sigma_{\alpha\in \Delta_j^{\theta, \delta}}\; 
\nu_{\alpha}(j)\,\hat{e}_{\alpha}(j)
\geq \int_M \tilde{\rho}^{\theta,\delta}(j) \omega (j)^n - 
\frac{C_6}{m(j)^{(1/2)-(3/2)\varepsilon}}.
\leqno{(4.13)}
$$
We next consider the function $\rho^{\theta,-}(j)$ and integrate it over $M$.
Then
\begin{align*}
&\rho^{\theta,-}(j) - \tilde{\rho}^{\theta,-}(j) \\
&\geq \;\frac{-S^{\theta,-}(j)Z^{\theta,+}(j)+\{Y^{\theta,-}(j)R^{\theta,-}-S^{\theta,-}(j)Z^{\theta,-}(j)\}}{Z(j)R(j)}.
\end{align*}
Put $f(t) := \log (\Sigma^{}_{\alpha \in\Delta^{\theta,-}}
|\tau'_{\alpha}|^2_{h(j)} \exp \{-t\, \epsilon_{\alpha}(j)\})$, 
and its second derivative $\ddot{f}(t)$ with respect to $t$
is nonnegative, wherever defined, 
and hence $\dot{f}(-2s(j)) \geq \dot{f}(0)$, i.e., 
$Y^{\theta,-}(j)R^{\theta,-}-S^{\theta,-}(j)Z^{\theta,-}(j) \geq 0$.
Thus
$$
\rho^{\theta,-}(j) - \tilde{\rho}^{\theta,-}(j) 
\; \geq \; \frac{-S^{\theta,-}(j)Z^{\theta,+}(j)}{Z(j)R(j)}.
$$
Integrating this over $M$, by (4.6) and (4.12), we obtain
$$
\int_M \rho^{\theta,-}(j) \omega (j)^n \;\geq\;  
\int_M \tilde{\rho}^{\theta,-}(j) \omega (j)^n\, 
- \,\frac{C_7}{m(j)^{(1/2)-\varepsilon}}.
\leqno{(4.14)}
$$
We now put ${a}_j (\theta ) := \int_M {\rho}^{\theta,-}(j)\omega (j)^n$ and
$\tilde{a}_j (\theta ) := \int_M \tilde{\rho}^{\theta,-}(j)\omega (j)^n$.
Then by setting
$$
a_{\infty} (\theta ) := \varlimsup_{j \to \infty} a_j ( \theta )
\quad \text{ and }\quad 
\tilde{a}_{\infty} (\theta ) := \varlimsup_{j \to \infty} \tilde{a}_j ( \theta ),
$$
we obtain from (4.14) the inequalities
 $$
 a_{\infty} (\theta ) 
\; \geq \; \tilde{a}_{\infty} (\theta )\; \geq \; 0
\leqno{(4.15)}
$$ 
for all $\theta$.
Since ${a}_j (0 ) = - \Sigma_{\alpha \in \Delta_j^-}\,
\nu_{\alpha}\hat{e}_{\alpha}(j) = r_j$, it now follows from (4.3) and (4.15) 
that
$$
r_{\infty} \; =\;  a_{\infty}(0)\; \geq\; \tilde{a}_{\infty}(0).
$$
Hence if  $\tilde{a}_{\infty}(0) > 0 $, then $r_{\infty} > 0$.
It now suffices to consider the remaining case $\tilde{a}_{\infty} (0) =0$.
Consider the positive constants $C_9$, $C_{10}$ $C_{11}$ in Lemma 5.1.
In view of $C_9 < C_{10}$, we put $C_{8}:= C_{10} - C_9$ for simplicity.
For every $a\in \Bbb R$, by using the Gauss symbol, let $[a ]$ denote the largest integer which does not 
exceed $a$.  Put $\gamma_j:= [C_{8} /\delta ]$.
Then
$$
C_{8}\log m(j)-1 < \gamma_j \leq C_{8} \log m(j).
$$
Now by setting $\theta_{\ell} := C_9 + (\ell - 1) \delta $, 
we obtain
$$
C_9 \leq
\theta_{\ell} < \theta_{\ell} + \delta 
\leq C_{10},   \qquad \ell  = 1,2,\dots, \gamma_j.
$$
Hence by applying Lemma 5.1 to $\theta = \theta_{\ell}$, we see from (4.13) 
and (5.2) the following inequality:
$$
\Sigma_{\alpha\in \Delta_j^{\theta_{\ell}, \delta}}\;
 \nu_{\alpha}(j)\,\hat{e}_{\alpha}(j)
\;\geq\; C_{11}\, \theta_{\ell}\, \delta \, - \,
\frac{C_6}{m(j)^{(1/2)-(3/2)\varepsilon}},
\quad \ell  = 1,2,\dots, \gamma_j.
$$
By summing these up, 
we obtain
\begin{align*}
r_j \,&=\, \Sigma^{}_{\alpha \in \Delta_j^+}\nu_{\alpha}(j)\, \hat{e}_{\alpha}(j) 
\,
\geq \;\Sigma_{\ell =1}^{\gamma_j}\,
\Sigma_{\alpha \in \Delta_j^{\theta_{\ell},\delta}}\; 
\nu_{\alpha}(j)\,\hat{e}_{\alpha}(j)\\
&\geq \, C_{11} \,\delta \;\Sigma_{\ell =1}^{\gamma_j}\, \theta_{\ell} 
\, -\, \frac{C_6 \,\gamma_j}{m(j)^{(1/2)-(3/2)\varepsilon}}
\geq \, C_{9}C_{11}\, \delta\, \gamma_j \, -\, \frac{C_6 \,\gamma_j}{m(j)^{(1/2)-(3/2)\varepsilon}}\\
&\geq \; C_9C_{11} \frac{ C_8 \log m(j)-1}{\log m(j)} \; -\; \frac{C_6 \,C_8\log m(j)}{m(j)^{(1/2)-(3/2)\varepsilon}}
\end{align*}
Then letting $j\to \infty$, we now conclude that
$$
r_{\infty} = \varlimsup_{j\to \infty} r_j \geq\; C_8C_9C_{11} \; > \;0,
$$
as required.
\qed

\section{Appendix}
In this appendix, under the same assumption as in the preceding section (hence 
$\Psi_{1,X'}^{\operatorname{SL}}$ is nontrivial),
we shall show the following:

\medskip\noindent
{\bf Lemma 5.1.}  
{\em If $\tilde{a}_{\infty} (0) =0$, then there exist positive real constants
$C_9$, $C_{10}$, $C_{11}$ independent of $j$ satisfying $C_9 < C_{10}$ such that
$$
\int_M \tilde{\rho}^{\theta,\delta}(j) \omega (j)^n \; \geq \;  C_{11} \,\theta\,\delta
\leqno{(5.2)}
$$
for all $\theta\in \Bbb R$ with $C_9\leq \theta < \theta + \delta \leq C_{10}$.
}

\medskip\noindent
{\em Proof}.
As in (3.15), we choose a basis $\{\tau_1, \tau_2, \dots, \tau_{N_1}\}$ for 
$(E_1)_0 = H^0(\mathcal{M}_0, \mathcal{L}_0)$ such that
$$
X'\cdot \tau_{\alpha} \; =\; e'_{\alpha}\, \tau_{\alpha},
\qquad \alpha = 1,2,\dots, N_1,
$$
where $e'_{\alpha}$ are real numbers
satisfying $\Sigma_{\alpha = 1}^{N_1}e'_{\alpha} = 0$. 
Since $\dim_{\Bbb C} \mathcal{M}_0 = n$,
the very ampleness of $\mathcal{L}_0$ implies that
$N_1 \geq n+1$.
Then we may assume without loss of generality that 
$e'_1 \leq e'_2  \leq \dots  \leq e'_{n+1}$ 
and that
$\tau_1$, $\tau_2$, \dots, $\tau_{n+1}$ are transcendental 
over $\Bbb C$ in the graded algebra
$$
\bigoplus_{m=1}^{\infty}\; H^0(\mathcal{M}^{}_0,\mathcal{L}_0^{\,m}).
$$
For later purposes,  consider the first positive integer $k$ satisfying the inequality 
$e_k < e_{k+1}$. 
Note that $k \leq n$.
Since $\Psi_{1,X'}^{\operatorname{SL}}$ is nontrivial,
there is an irreducible component $\mathcal{F}$ of $\mathcal{M}_0$
with subscheme structure in $\mathcal{M}_0$
such that $T'$ acts on the scheme $\mathcal{F}$ nontrivially.
Hence, in view of the fact that the generically finite rational
 map $\tau : \mathcal{F} \to \Bbb P^n(\Bbb C )$ defined by
$$
\tau (p) := (\tau_1 (p): \tau_2 (p): \dots : \tau_{n+1}(p)),
\qquad p \in \mathcal{F},
$$
is $T'$-invariant,  
the set $\{e'_{\alpha}\,; \,\alpha = 1,2,\dots,n+1\}$ 
contains at least two distinct real numbers.
In particular, 
$$
e'_1 \; <\; e'_{n+1}.
$$
Let $B(j)$ be the of all $b =(b_1, b_2, \dots, b_{n+1})
$$\in \Bbb Z^{\,n+1}$ 
such that
$$
b_1 \geq 0, \; b_2\geq 0,\; \dots,\; b_{n+1}\geq 0,\; 
b_1 +b_2 + \dots +b_{n+1} = m(j).
$$
Then the weight $w(b)$
of the action of $X'$ on 
$\tau^b := \tau_1^{b_1}\tau_2^{b_2}\dots \tau_{n+1}^{b_{n+1}}$ is 
$$
w(b)\; =\; k_j\, +\,
\Sigma_{\alpha =1}^{n+1}\; b_{\alpha}e'_{\alpha}, \qquad b\in B(j),
$$
where $k_j$ is a real constant 
independent of the choice of $b$
such that the sequence $k_j/m(j)$, $j$=1,2,\dots, is bounded.
Here $k_j$ is chosen in such a way that the representation matrix of 
$X'$ on $(E_{m(j)})_0$ 
is traceless.
In view of (3.16), replacing
$$
c (j):= \; \frac{k_{j}}{m(j)} \,-\, \bar{\xi}(j), \qquad j=1,2,\dots,
$$ 
by its subsequence if necessary,
we may assume without loss of generality that $c(j)$
 converges to some real number $c (\infty )$
as $ j\to \infty$. For each $\theta \in \Bbb R$ with $0 \leq \theta  \leq C_2$,
we define $\epsilon_b (j)_{\theta}$ and $\hat{\epsilon}_b (j)_{\theta}$ by 
$$
\begin{cases}
\;\;\epsilon_b (j)_{\theta} &:= w(b) - m(j) \bar{\xi} (j) - m(j) \theta \\
&\,= m(j)\{ c(j)\, -\, \theta\}
 + \, \Sigma_{\alpha =1}^{n+1}\; b_{\alpha}e'_{\alpha},\\
\;\;\hat{\epsilon}_b (j)_{\theta} &:= \epsilon_b (j)_{\theta}/m(j).
\end{cases}
$$
Note that ${\epsilon}^{\operatorname{max}}(j)^{}_{\theta} 
:= \max^{}_{b\in B(j)} \epsilon_b (j)_{\theta}$ 
is exactly
$m(j) \{ c (j)  +  e'_{n+1} - \theta \}$.
Put $\eta := c (\infty ) + e'_{n+1}$. Then one of the following 
two cases occurs:

\medskip
\quad Case 1: $\eta \leq 0$, 
\qquad Case 2: $\eta >0$.

\medskip\noindent

\noindent
First, we consider Case 1, 
i.e., let $\eta \leq 0$. 
 Let $\varepsilon_1 > 0$ be a sufficiently small real 
number independent of the choice of $j$. Then 
\begin{align*}
{\epsilon}^{\operatorname{max}}(j)^{}_0  - \varepsilon_1\, m(j)\;
&=\;  \{\,\eta - c(\infty ) +
c(j)  - \varepsilon_1\,\} \, m(j) \tag{5.3}
\\
&\leq \; \{- c(\infty ) +
c(j)  - \varepsilon_1\,\} \, m(j) 
\leq \; 0
\end{align*}
for $j \gg 1$. 
Since $\int_M |\tau'_{\alpha}|^2_{h(j)} \omega (j)^n$ is $1$ for all $\alpha$ and $j$,
then from
$$
\tilde{a}_j (0) = \int_M \tilde{\rho}^{0,-}(j) \,\omega (j)^n
\; =\; \int_M \frac{\Sigma_{\alpha\in \Delta_j^-}\, \{-\hat{e}_{\alpha}(j)\}|\tau'_{\alpha}|^2_{h(j)}}{\Sigma_{\alpha =1}^{N(j)}|\tau'_{\alpha}|^2_{h(j)}}\, \omega (j)^n,
$$
and also from 
$$
\Sigma_{\alpha =1}^{N(j)}|\tau'_{\alpha}|^2_{h(j)} 
= C_{12} m(j)^n\{1+ O(1/m(j))\},
\leqno{(5.4)}
$$
it now follows that
$$
\tilde{a}_j (0)\; =\;  \frac{C_{13}\left \{1+ O(1/m(j)) \right \} \Sigma^{}_{\alpha\in \Delta_j^-}
\,\{-{e}_{\alpha}(j)\}}{m(j)^{n+1}}. 
\leqno{(5.5)}
$$
Let $\Sigma_{-}\{-\epsilon_b (j)_0\}$ denote the sum of all $-\epsilon_b (j)_0$ 
with $b$ running through the set of all  $b$ in $B(j)$ 
such that $\epsilon_b (j)_0 < 0$.
Then
\begin{align*}
&\Sigma_{\alpha\in \Delta_j^-} 
\,\{-{e}_{\alpha}(j)\}\; \geq \; \Sigma_-\,\{-\epsilon_b (j)_0\} \;
\geq \; \Sigma_{b\in B(j)}\, \{- \epsilon_b (j)_0\}  \tag{5.6}
\\
&\geq \; \Sigma_{b\in B(j)}\, \{\epsilon^{\operatorname{max}} (j)_0- \epsilon_b (j)_0\}
\,-\, \Sigma_{b\in B(j)}\, \epsilon^{\operatorname{max}} (j)_0\\
& \geq \; \Sigma_{b\in B(j)}\, \{\epsilon^{\operatorname{max}} (j)_0- \epsilon_b (j)_0\}
\,-\, \Sigma_{b\in B(j)}\, \varepsilon_1 m(j),
\end{align*}
where the last inequality follows from (5.3).
The first term  in the last line is estimated as follows:
\begin{align*}
&\Sigma_{b\in B(j)}\, \{\epsilon^{\operatorname{max}} (j)_0- \epsilon_b (j)_0\}
\tag{5.7} \\
&\geq\; \Sigma_{b\in B(j)} 
\Sigma_{\alpha = 1}^{n+1}\, b_{\alpha} (e'_{n+1}- e'_{\alpha})\;
\geq  \; (e_{n+1}'-e'_1)\,\Sigma_{b\in B(j)}\; b_1 \\
&=\; \frac{e'_{n+1}- e'_1}{(n-1)!}\,\Sigma_{b_1 =0}^{m(j)} \left \{  
b_1\,\Pi_{\ell =1}^{n-1}\,( m(j)-b_1+\ell ) \right \}\\
&= \; \frac{e'_{n+1} - e'_1}{(n-1)!}\;\Sigma_{b_1 =0}^{m(j)}\; b_1\left \{
(m(j) -b_1)^{n-1} + O(m(j)^{n-2}) \right \} \\
&=\;\frac{e'_{n+1}-e'_1}{(n+1)!}\, m(j)^{n+1}
+ O(m(j)^n).
\end{align*}
Since 
$\,\Sigma_{b\in B(j)}\, \varepsilon_1\, m(j) \; =\; (\varepsilon_1/n!)\, 
\Pi_{\ell =0}^{n}\; (m(j) + \ell )$, it now follows from (5.5), (5.6) 
and (5.7) that
$$
\tilde{a}_j(0) \; \geq  \; C_{13}\left \{ \frac{e'_{n+1} - e'_1}{(n+1)!} - 
\frac{\varepsilon_1}{n!} \right \} \,  + O (1/m(j)).
$$
Since $\varepsilon_1$ is sufficiently small, we obtain in this case
$$
\tilde{a}_{\infty}(0)\; =\; \varlimsup_{j\to \infty} \tilde{a}_j (0)\; 
\geq \; C_{13}\left \{ \frac{e'_{n+1} - e'_1}{(n+1)!} - 
\frac{\varepsilon_1}{n!} \right \}  \; >\; 0,
$$
in contradiction to the assumption $\tilde{a}_{\infty}(0) = 0$.

\medskip
Hence, it suffices to consider Case 2. 
In this case by $\eta > 0$, using  the Gauss symbol, we consider the nonnegative integer 
$p:= [\,(e_{n+1}' - e_1')/\eta\,]$. 
We then choose positive constants $C_9$ and $C_{10}$ by 
$$
C_9 := \eta -\frac{e_{n+1}' - e_1'}{p+3}
\quad\text{ and }\quad C_{10} := \eta -\frac{e_{n+1}' - e_1'}{p+4}.
$$
Take an arbitrary real number $\theta$ such that $C_9 \leq \theta < \theta +\delta \leq C_{10}$. 
Then for some $p^{}_{\theta}\in \Bbb R$ satisfying 
$p+3 \leq  p^{}_{\theta} < p+4$, we can write $\theta$ as 
$$
\theta \; =\; \eta -\frac{e_{n+1}' - e_1'}{p^{}_{\theta}}
$$
Again by using the Gauss symbol, we put
\begin{align*}
v^{}_{j,\theta}\, &:= \frac{c(j)-c(\infty )}{e_{n+1}' - e_1'} + \frac{1}{p^{}_{\theta}},
\tag{5.8} \\
m' \; &:=\; [\,m(j)\, v^{}_{j,\theta}\, ].
\tag{5.9}
\end{align*}
Here $1/(p+5) \leq  v_{j,\theta}\leq 2/5$
for all $j \geq j_0$,  where $j_0$ is a sufficiently large positive integer 
independent of the choice of $\theta$.
Put $C_{14}:=(e'_{n+1} - e'_1)/5$ for simplicity.
For $j\gg 1$, let $B'(j)$ denote the subset of $B(j)$ consisting of all $b = (b_1,b_2, \dots, b_{n+1}) \in B(j)$ such that $b_1$, $b_2$, \dots, $b_{n+1}$ are expressible as
$$
\begin{cases}
&b_1  = m' + (1/5)m(j) - \Sigma_{\alpha =2}^{n+1} \, \beta_{\alpha},\\
&b_{\alpha} = \beta_{\alpha} ,\qquad \alpha =2,3,\dots, n,\\
&b_{n+1} = (4/5)m(j) - m' + \beta_{n+1},
\end{cases}
$$
for some $\beta_2$, $\beta_2$, \dots, $\beta_{n+1}$, 
where $\beta = (\beta_2, \beta_3,\dots, \beta_{n+1})$  runs through the 
subset $\,B''(j)\,$ of $\Bbb Z^{\,n}_{\geq 0}$ consisting of all $\beta\in \Bbb Z^n$ such that
$$
\begin{cases}
&\beta_2 \geq 0,\; \beta_3\geq 0,\; \dots,\;
\beta_{n+1}\geq 0,\;\Sigma_{\alpha =2}^{k} \, \beta_{\alpha} \leq  (1/5)m(j) ,
\\
&C_{14}\, m(j) \leq\Sigma_{\alpha =k+1}^{n+1}\, \beta_{\alpha}(e_{\alpha}'-e'_1)
\leq  (C_{14} +\delta )\, m(j) - (e_{n+1}'-e_1'), 
\end{cases}
$$
where for $k=1$, the condition 
$\Sigma_{\alpha =2}^{k} \, \beta_{\alpha} \leq  (1/5)m(j)$ is assumed to be void. 
Then for each $b = (b_1, b_2, \dots , b_{n+1})\in B'(j)$, 
\begin{align*}
\epsilon^{}_{b}(j)^{}_{\theta}
\, &=\,\epsilon^{\operatorname{max}}(j)^{}_{\theta}\,- \left \{ \epsilon^{\operatorname{max}}(j)^{}_{\theta}- \epsilon^{}_{b}(j)^{}_{\theta}\right \} 
 \\
&= \,  m(j)\{\eta - c(\infty ) + c(j) - \theta \}\,-\,\Sigma_{\alpha =1}^{n+1}\, b_{\alpha} (e'_{n+1} - e'_{\alpha}) 
\\
&= \, m(j)\, v^{}_{j,\theta}\,(e_{n+1}'-e_1')\,
-\,\Sigma_{\alpha =1}^{n+1}\, b_{\alpha} (e'_{n+1} - e'_{\alpha})\\
&=\, \{ m(j) \, v^{}_{j,\theta}  -m' -(1/5)m(j)\}(e_{n+1}' - e_1')\, +\,
 \Sigma_{\alpha = 2}^{n+1} \,\beta_{\alpha} (e_{\alpha}'-e_{1}') \\
&\geq\,  -C_{14}m(j) + \Sigma_{\alpha = k+1}^{n+1} \,\beta_{\alpha} (e_{\alpha}'-e_{1}')
\, \geq \, 0,
\end{align*}
where we used the inequality $m(j) v^{}_{j,\theta}  -m' \geq 0$ (cf. (5.9)).
Again by  (5.9), $m(j) v^{}_{j,\theta}  -m' <1$ and hence, for each $b \in B'(j)$
\begin{align*}
\epsilon^{}_{b}(j)^{}_{\theta}\; &< \;\{1 -(1/5)m(j)\}\,( e_{n+1}' - e_1')\, + \,
\Sigma_{\alpha = 2}^{n+1} \,\beta_{\alpha} (e_{\alpha}'-e_{1}')\\
&\leq \;e_{n+1}' - e_1' - C_{14} m(j) + \Sigma_{\alpha =k+1}^{n+1}\, 
\beta_{\alpha} (e_{\alpha}'-e_{1}')\;
\leq\; m(j)\delta.
\end{align*}
Hence by (5.4) and $-C_0(\log m(j)) q(j) \leq s(j) \leq 0$, it follows that 
\begin{align*}
&\int_M \tilde{\rho}^{\theta,\delta}(j)\omega (j)^n
\;=\; C_{15}\,\theta \delta^{-1}\frac{\{ 1+ O(1/m(j))\}}{m(j)^n}
\int_M \tilde{Y}^{\theta,\delta} \omega (j)^n\\
&=\;  C_{15}\,\theta \delta^{-1}\frac{\{ 1+ O(1/m(j))\}}{m(j)^n}
\Sigma^{}_{\alpha \in \Delta_j^{\theta,\delta}}\,\hat{\epsilon}_{\alpha}(j)
\exp\{2s(j)\epsilon_{\alpha}(j)\}\\
& \geq \; C_{15}\,\theta \delta^{-1}\frac{\{ 1+ O(1/m(j))\}}{m(j)^n}
\Sigma^{}_{b \in B'(j)}\,\hat{\epsilon}_{b}(j)_{\theta}\,
\exp\{2s(j)\epsilon_{b}(j)_{\theta}\}\\
&\geq \; C_{16}\,\theta \delta^{-1}\frac{\{ 1+ O(1/m(j))\}}{m(j)^n}
\Sigma^{}_{b \in B'(j)}\,\hat{\epsilon}_{b}(j)_{\theta},
\end{align*}
where the last inequality follows from $2s(j){\epsilon}_{b}(j)_{\theta}\geq 
- 2 C_0$. Put 
$$
\hat{x}(\beta ):= \frac{
\Sigma_{\alpha = k+1}^{n+1} \,\beta_{\alpha} (e_{\alpha}'-e_{1}') 
- C_{14}m(j)}{m(j)}
$$ 
for simplicity.
Now by $\epsilon^{}_{b}(j)^{}_{\theta} \geq  \hat{x}(\beta ) m(j)$, we see that
$$
\int_M \tilde{\rho}^{\theta,\delta}(j)\omega (j)^n
\; \geq\; C_{16}\,\theta \delta^{-1}\frac{\{ 1+ O(1/m(j))\}}{m(j)^n}
\; \Sigma_{\beta\in B''(j)} \hat{x}(\beta ).
\leqno{(5.10)}
$$
For $j \to \infty$, as far as the growth order of $\Sigma_{\beta\in B''(j)} \hat{x}(\beta )$ 
with respect to $m(j)$ is concerned, we may assume without loss of generality
$$
e'_{k+1}- e_1' \;=\; e_{k+2}'-e_1' \;=\; \dots \; =\; e_{n+1}'-e_1'
$$ 
Then a computation similar to (5.7) 
allows us to write $\Sigma_{\beta\in B''(j)} \hat{x}(\beta ) $ 
as
\begin{align*}
&C_{17} m(j)^{k-1}\cdot 
\{\delta m(j)\}^2 m(j)^{n-k-1} + \text{lower order term in $m(j)$}\\
& = \; C_{17}\,\delta^2  \,m(j)^{n}+ \text{lower order term in $m(j)$}.
\end{align*}
This together with (5.10) implies (5.2), as required.
\qed


\bigskip\noindent
{\footnotesize
{\sc Department of Mathematics}\newline
{\sc Osaka University} \newline
{\sc Toyonaka, Osaka, 560-0043}\newline
{\sc Japan}}
\end{document}